\documentclass[notitlepage,leqno,11pt]{article}
\usepackage{amssymb}
\catcode`\@=11 \@addtoreset{equation}{section}

\catcode`\@=12

\usepackage{latexsym}
\usepackage{textcomp}
\usepackage{amsmath}
\usepackage{amsfonts}
\usepackage{amssymb}
\usepackage{mathrsfs}

\renewcommand{\a }{\alpha }
\renewcommand{\b }{\beta }
\renewcommand{\d}{\delta }

\newcommand{\D }{\Delta }

\newcommand{\e }{\varepsilon }
\newcommand{\g }{\gamma}

\newcommand{\G }{\Gamma}
\renewcommand{\l }{\lambda }

\newcommand{\n }{\nabla }
\newcommand{\vp }{\varphi }
\newcommand{\vt }{\vartheta }

\newcommand{\s }{\sigma }

\renewcommand{\t }{\tau }
\newcommand{\z }{\zeta}

\renewcommand{\O }{\Omega }

\newcommand{\ov}{\overline}

\newcommand{\be}{\begin{equation}}
\newcommand{\ee}{\end{equation}}

\newcommand{\R}{\mathbb{R}}
\newcommand{\N}{\mathbb{N}}

\newcommand{\de}{\partial}

\newcommand{\ti}{\widetilde}

\newcommand{\ra}{{\rangle}}
\newcommand{\la}{{\langle}}

\renewcommand{\k}{\kappa}

\newcommand{\scrD }{\mathscr{D}}

\newcommand{\scrH }{\mathscr{H}}

\newcommand{\calH }{\mathcal{H}}

\newcommand{\calL }{\mathcal{L}}

\newcommand{\calD }{\mathcal{D}}

\newcommand{\calB }{\mathcal{B}}

\newcommand{\calF}{{\mathcal F}}

\newcommand{\calS}{{\mathcal S}}

\newtheorem{Theorem}{Theorem}[section]
\newtheorem{Lemma}[Theorem]{Lemma}
\newtheorem{Proposition}[Theorem]{Proposition}

\newtheorem{Remark}[Theorem]{Remark}
\newtheorem{Definition}[Theorem]{Definition}

\def\proof{\noindent{{\bf Proof. }}}
\def\square{\vbox{
    \hrule height .4pt
    \hbox{\vrule width .4pt height 7pt \kern 7pt
       \vrule width .4pt}
    \hrule height .4pt }}

\def\square{\vbox{
    \hrule height .4pt
    \hbox{\vrule width .4pt height 7pt \kern 7pt
       \vrule width .4pt}
    \hrule height .4pt }}

\def\QED{\hfill {$\square$}\goodbreak \medskip}

\def\R{{\mathbb R}}

\def\f{{\varphi}}
\def\div{{\rm div}}

\newcommand{\Ds}{(-\Delta)^s}

\font\sc=cmcsc9  \textwidth=14truecm
\begin{document}

\title{ Semilinear elliptic equations for the  fractional Laplacian  with
Hardy potential}

\author{Mouhamed Moustapha Fall\footnote
{\footnotesize{Goethe-Universit\"{a}t
Frankfurt, Institut f\"{u}r Mathematik.  Robert-Mayer-Str. 10, D-60054 Frankfurt am Main, Germany. 
 E-mail: {\tt fall@math.uni-frankfurt.de,  mouhamed.m.fall@gmail.com}. } } }
\date{}

\maketitle

\bigskip

\noindent {\footnotesize{\bf Abstract.} In this paper we study existence and  nonexistence of
nonnegative distributional solutions for a class of
semilinear fractional elliptic equations  involving  the Hardy  potential.}
\bigskip\bigskip



\section*{Introduction}
Let $B$ be a ball of   $\R^N$, $N> 2s$, centered at 0. Let $s\in(0,1)$, $p>1$ and $\g\geq0$. In this paper, we study
existence and nonexistence of nonnegative functions $u\in\calL^1_s\cap L^p_{loc}(B)$ satisfying
\be\label{eq:inequalityb}
\Ds u- \g |x|^{-2s} u =u^p\quad\textrm{ in } B,
\ee
where
$$
\calL^1_s=\left\{u:\R^N\to\R\,:\, \int_{\R^N}\frac{|u|}{1+|x|^{N+2s}}<\infty\right\}.
$$
Equality \eqref{eq:inequalityb} is understood in the sense of distributions. The distribution $\Ds u $ is defined as
$$
\la\Ds u,\vp \ra =\int_{\R^N} u\Ds\vp dx \quad\forall \vp \in C^\infty_c(B).
$$
Here the fractional Laplacian  $\Ds$  is defined via the Fourier transform as
\be\label{eq:def-frac_lap_four-int}
{\Ds\vp}(x) = \frac{1}{(2\pi)^{\frac{N}{2}}}\,
\int_{\R^N} |\zeta|^{2s}\widehat{\vp}(\zeta) e^{\imath\zeta\cdot x }d\zeta,
\ee
where
$$
\widehat{\vp}(\zeta)=\calF(\vp)(\zeta)=\frac{1}{(2\pi)^{\frac{N}{2}}}\,
\int_{\R^N}e^{-\imath\zeta\cdot{x} }\vp(x)dx
$$
is the Fourier transform of $\vp$. The nonlocal structure of  the  fractional Laplacian  $\Ds$
can be seen in its representation in the real space:
\be\label{eq:DSvp-int}
\Ds \vp(x)=C_{s,N}\,P.V.\int_{\R^N}\frac{\vp(x)-\vp(y)}{|x-y|^{N+2s}}dy,
\ee
for some positive  constant $C_{s,N}$.
For the equivalence between \eqref{eq:DSvp-int} and \eqref{eq:def-frac_lap_four-int}, we refer the reader to \cite{Landk}.\\
Problem \eqref{eq:inequalityb} is related to the relativistic
 Hardy inequality which were proved by Herbst in \cite{Herb} (see also \cite{Yaf}):
\be\label{eq:fracHardy}
\g_0\int_{\R^N}|x|^{-2s} u^2dx\leq \int_{\R^N}|\z|^{2s} \widehat{u}^2d\z  \quad \forall u\in C^\infty_c(\R^N),
\ee
where
\be\label{eq:defgamzero}
\g_0=2^{2s} \frac{\Gamma^2\left(\frac{N+2s}{4}\right)}{\Gamma^2\left(\frac{N-2s}{4}\right)}.
\ee
%
%
The constant $\g_0$ is optimal and converges to  the classical Hardy constant $\frac{(N-2)^2}{4}$ when $s\to1$.
Here $\Gamma$ is the usual gamma function. We should mention that \eqref{eq:fracHardy} is a particular case of the Stein and Weiss inequality, see \cite{StWe}.\\

A great deal of work is currently been  devoted to the study of the fractional Laplacian  as it appears in  many fields
such as probability theory, physics and mathematical finance. We refer the reader  to papers
 \cite{CS}, \cite{Sil}, \cite{SV}, \cite{FQT}, \cite{BBC} (and the references there in)
 for a nice expository. A good reference for the  potential theory of $\Ds$ can be found in the book of Landkof \cite{Landk}.
 The operator $ \Ds-\g_0|x|^{-2s}$ appears in the problem
of stability of relativistic matter in magnetic fields. One can see  \cite{Rup}  where
a lower bound and a Gagliardo-Nirenberg-type  inequality were  proved.\\
The problem of existence and nonexistence of  \eqref{eq:inequalityb}, for $s=1$,
 was studied by Brezis-Dupaigne-Tesei in \cite{BDT}
where the authors showed that for $\b\in\left[0,\frac{N-2}{2}\right]$, $ 1<p<\frac{N+2-2\b}{N-2-2\b}$, the problem
$$
-\D u-\left(\frac{(N-2)^2}{4}-\b^2\right)|x|^{-2}=u^p\quad\textrm{ in } \calD'(B)
$$
has a positive solution $u\in L^p(B)$ and does not have any nonnegative and nontrivial  supersolution
$u\in L^p_{loc}(B\setminus\{0\})$
when $\b\in\left[0,\frac{N-2}{2}\right)$ and  $p\geq\frac{N+2-2\b}{N-2-2\b}$. Some related results and problems
are in \cite{BC}, \cite{BDT}, \cite{DD-H}, \cite{D}, \cite{DN}, \cite{Fall-ne-sl}, \cite{T}, \cite{FaMu1}, \cite{FaMu-ne}, \cite{mmf}, \cite{BM}, \cite{BMS}.\\
Our results in this paper extends the one  of
 Brezis-Dupaigne-Tesei in \cite{BDT} to the case $s\in(0,1)$. Before stating them, we fix the following notation:
for  $\a\in\left[0,\frac{N-2s}{2}\right)$, we put
$$
\g_\a=2^{2s}  \frac{\Gamma\left(\frac{N+2s+2\a}{4}\right)}{\Gamma\left(\frac{N-2s-2\a}{4}\right)} 
\frac{\Gamma\left(\frac{N+2s-2\a}{4}\right)}{\Gamma\left(\frac{N-2s+2\a}{4}\right)}.
$$
The mapping   $\a\mapsto\g_\a$ is monotone decreasing and $\g_\a\to0$ when $\a\to \frac{N-2s}{2}$.
 We should mention that  this constant appears in the perturbation of the
 'ground-state' $|x|^{\frac{2s-N}{2}}$ for the operator $ \Ds-\g_0|x|^{-2s}$. Indeed, letting $\vartheta_\a=|x|^{\frac{2s-N}{2}+\a}$
we have
$$
 \Ds\vartheta_\a-\g_\a|x|^{-2s}\vartheta_\a=0\quad\textrm{ in } \R^N\setminus\{0\},
$$
see Lemma \ref{lem:grstate} in Section \ref{s:nonexist}.\\
Our existence result is   the following
\begin{Theorem}\label{thm:main-th-ex}
Let $\a\in\left[0,\frac{N-2s}{2}\right]$ and $1<p< \frac{N+2s-2\a}{N-2s-2\a}$. There exits a function $u \in \calL^1_s\cap L^p(B)$ satisfying
$ u>0$ in $B$ and
$$
(-\Delta)^s u-\g_\a |x|^{-2s}\, u = u^p\quad\textrm{ in } \calD'(B).
$$
\end{Theorem}
As what concerns nonexistence, we have obtained:
\begin{Theorem}\label{thm:main-th-ne}
Let  $\a\in\left[0,\frac{N-2s}{2}\right)$ and  $u \in \calL^1_s\cap L^p_{loc}(B\setminus\{0\})$ such that  $u\geq 0$  and
$$
(-\Delta)^s u-\g_\a |x|^{-2s}\, u \geq u^p \quad\textrm{ in }\calD'(B\setminus\{0\}).
$$
If $p\geq  \frac{N+2s-2\a}{N-2s-2\a} $, then $u=0$ in $B$.
\end{Theorem}

 We observe that if  $\a=0$ we have $p+1=\frac{2N}{N-2s}$:
the critical Hardy-Littlewood-Sobolev exponent and that $\frac{N+2s-2\a}{N-2s-2\a}\to+ \infty $ as $\a\to \frac{N-2s}{2}$.\\
 The proof of Theorem \ref{thm:main-th-ne} relies
on weak comparison principles recently used by the author in \cite{Fall-ne-sl}.  However
 substantial difficulties have to be overcome due to the nonlocal structure of the fractional Laplacian.
Nonexistence result of nonlinear elliptic  problems using comparison principles have been obtained in
\cite{AS-CPDE}, \cite{AS-RIMS} \cite{PoTe},
\cite{LLM-Proce}, \cite{KLS},
\cite{KLS-JDE}, \cite{KLS-Trans} and the references therein.\\
For the  existence result, in the supercritical case $\frac{N+2s-2\a}{N-2s-2\a}>p\geq \frac{N+2s}{N-2s}$,
we have an explicit solution constructed via $\vartheta_\a$. In  the subcritical case,  $p+1<\frac{2N}{N-2s}$,
we used standard  variational arguments thanks to the following improved fractional Hardy  inequality:
\begin{Theorem}\label{th:improv-frHa}
 Let  $ 2>q>\max\left(1,\frac{2}{1+2s}\right) $.
Then there exists a constant $C>0$ such that for all $u\in C^\infty_c(B)$,
\be
 C \|u\|^{2}_{W^{\t,q}_0(B )} \leq \g_0 \int_{B}|x|^{-2s}u^2dx- \int_{\R^N}|\z|^{2s}\widehat{u}^2d\z  ,
\ee
where  $ \t=\frac{1+2s}{2}-\frac{1}{q}$.
\end{Theorem}
This result, which might be of self interest, is proved in  Appendix \ref{a:remainder}.\\
The proof of all the results  presented  above are manly  based on a
 Dirichlet-to Neumann operator $\calB_s$  for which $\calB_s u =\Ds \ti{u}$ in $\calD'(E)$ for any Lipschitz bounded open set $E$ of $\R^N$
and $\ti{u}$ is the null extension outside $E$ of a function  $u$ belonging to  some Sobolev space.
To be more precise let us first recall  the result of Caffarelli and Silvestre.
We recall that
$$
H^s(\R^N)=\left\{u:\R^N\to\R \,:\, (1+ |\zeta|^s)\widehat{u}\in   L^2(\R^ N)\right\}.
 $$
Put
$$
\R^{N+1}_+=\{(t,x)\,: \, t>0,\, x\in \R^N\}.
$$
Given $w\in H^s(\R^N)$, minimization procedure yields the existence of a unique function  $\calH(w)\in H^1(\R^{N+1}_+; t^{1-2s} )$
 being   the harmonic extension of $w$ over the half space $\R^{N+1}_+$:
\be\label{eq:Harmext-int}
\begin{cases}
\div(t^{1-2s}\n \calH(w))=0\quad \textrm{ in }\R^{N+1}_+,\\
\calH(w) = w \quad \textrm{ on }\R^N.\\
\end{cases}
\ee
In  \cite{CSilv},   Caffarelli and Silvestre
proved that  $\Ds w$ is given by the  Dirichlet-to-Neumann operator $\lim_{t\to 0}t^{1-2s}\frac{\de \calH(w)}{\de t}$:
$$
-\lim_{t\to 0}t^{1-2s}\frac{\de \calH(w)}{\de t}= \k_{s} \Ds w \quad   \textrm{ in } \R^N,
$$
for some constant $\k_s>0$. In addition
\be\label{eq:ncaHeqns2-int}
\int_{\R^{N+1}_+}t^{1-2s}|\n \calH(w)|^2dxdt=\k_{s}\int_{\R^N}|\zeta|^{2s}\widehat{w}d\zeta.
\ee

 We want to provide similar arguments in bounded open sets.
We define  the Hilbert space $\scrD^{s,2}(\R^N)$ which is  the completion of  $C^\infty_c(\R^N)$
with respect to the norm:
$$
v\mapsto \int_{\R^N}|\zeta|^{2s}|\widehat{v}|^2 d\zeta.
$$
Let $E$ be a bounded  open set    in $\R^N$ with Lipschitz boundary.
We   introduce the Hilbert space
$$
\scrH^s_0(E):=\{ u\in H^s(E)\,:\,{\ti{u}} \in \scrD^{s,2}(\R^N) \},
$$
where
$$
\ti{u}=\begin{cases}
        u\quad \textrm{ in } E\\
        0\quad\textrm{ in } \R^N\setminus E.
       \end{cases}
$$
The space $\scrH^s_0(E)$ is endowed with  the natural norm
\be\label{eq:normscrH-int}
\| u \|_{\scrH^s_0(E) }^2=\int_{\R^N}|\zeta|^{2s}|\widehat{\ti{u}}|^2 d\zeta=\int_{\R^N}|\zeta|^{2s}|\calF(\ti{u})|^2 d\zeta.
\ee
Note that, since $E$ is bounded, by \eqref{eq:fracHardy}  there exists a constant $C(E)>0$ such that
\be\label{eq:Poincar-scrH}
C(E)\|\ti{u}\|_{ H^s(\R^N) }\leq \| u \|_{\scrH^s_0(E) }\leq \|\ti{u}\|_{ H^s(\R^N) }\quad \forall u \in \scrH^s_0(E).
\ee
From this we deduce that
\be\label{eq:scrh0eq}
\scrH^s_0(E)=\{ u\in H^s(E)\,:\, \ti{u} \in  H^s(\R^N)\}.
\ee
Hence, see for instance [\cite{Gris}, Theorem 1.4.2.2 ], the space $ C^\infty_c(E)$
is  dense in $\scrH^s_0(E)$.\\
By \eqref{eq:scrh0eq}  for any $u\in \scrH^s_0(E)$, we can consider its harmonic extension $\calH(\ti{u})$ as in \eqref{eq:Harmext-int}.
We define  the Dirichlet-to-Neumann operator $\calB_s:\scrH^s_0(E)\to  \scrH^{-s}(E) $ given by
$$
\calB_s u =-\k_s^{-1}\lim_{t\to 0}t^{1-2s}\frac{\de \calH(\ti{u})}{\de t},
$$
where $\scrH^{-s}(E) $ is the dual of $\scrH^s_0(E)$. This  operator turns out to be  linear and it is an   isometry,
$$
 \|\calB_s u\|_{\scrH^{-s}(E)}  =\|u\|_{\scrH^s_0(E) }\quad\forall u\in \scrH^s_0(E)
$$
by \eqref{eq:ncaHeqns2-int}.
Moreover
$$
\calB_s u= \Ds \ti{u}\quad \textrm{ in } \calD'(E)\quad \forall u\in \scrH^s_0(E).
$$
In particular a solution $u\in \scrH^s_0(E)$ to the  problem
$$
\calB_s u=  f\quad\textrm{ in } E
$$
yields a solution $v\in \scrD^{s,2}(\R^N)$ to  the problem
$$
\begin{cases}
  \Ds v= f\quad\textrm{ in } E,\\
  \hspace{1cm} v= 0 \,\,\quad\textrm{ in } \R^N\setminus E
\end{cases}
$$
and conversely. We refer to the next section for more details.\\
 In order to get, say, qualitative informations on the solution to the problem
$$
\calB_s u= \Ds \ti{u}= f,
$$
it is, in general,   more convenient to  work with the (mixed) problem
\be\label{eq:Harmext-intH}
\begin{cases}
\div(t^{1-2s}\n \calH(\ti{u}))=0\quad \textrm{ in }\R^{N+1}_+,\\
-\lim_{t\to 0}t^{1-2s}\frac{\de \calH(\ti{u})}{\de t}= \k_{s}f \quad   \textrm{ in }E.
\end{cases}
\ee
New difficulties arise here because of  the (possible) degeneracy of the  equation \eqref{eq:Harmext-intH}. However
the weight $ t^{1-2s}$ falls into the  Muckenhoupt class of weights thus regularity results, Harnack inequalities are
available (see \cite{FKS}) and this is enough for our purpose in this paper.\\
An interesting characterization of $\scrH^s_0(E) $, see  \cite{Gris},  is that $\scrH^s_0(E) $
is the interpolation space $(H^2_0(E),L^2(E))_{s,2}$:
\be\label{eq:cahrac-scrH}
\scrH^s_0(E)=\left\{ \begin{array}{c}
H^{s}(E)\quad s\in(0,1/2),\\
H^{s}_{00}(E)\quad s=1/2,\\
H^{s}_0(E)\quad s\in(1/2,1),
\end{array}\right.
\ee
where
$$
H^{\frac{1}{2}}_{00}(E)=\left\{u\in H^{\frac{1}{2}}(E)\,:\, \int_{E}\frac{u^2(x)}{d(x)}dx<\infty\right\},
$$
endowed with the natural norm,
with $d(x)=\textrm{dist}(x,\de E)$.
\begin{Remark}\label{rem:Cabre-Tan}
Let $E$ be a smooth bounded domain of $\R^N$. Recently a pseudo differential  operator   $A_s$
 of order $2s$ was introduced by Cabr\'e and Tan \cite{CT} for $s=1/2$  (see  \cite{CDDS} for every $s\neq1/2$)
 in the following way:
 for any $u\in \scrH^s_0(E)$
$$
A_s u=\sum_{k=1}^\infty \mu_k^s u_k\vp_k,
$$
where $\mu_k $ is the zero Dirichlet eigenvalues of $-\D$  with
 corresponding orthonormal eigenfunctions $\vp_k$ and $u_k=\int_{E}u\vp_kdx$
is the component of $u$ in the $L^2(E)$ basis $\{\vp_k\}$.\\
Using \eqref{eq:cahrac-scrH}, it was shown in \cite{CT} and \cite{CDDS} that
\be\label{eq:normCT}
\scrH^s_0(E)=\left\{ u\in L^2(E)\,:\, \sum_{k=1}^\infty \mu_k^s|u_k|^2<\infty\right\}.
\ee
The  operator $A_s$ corresponds to the Dirichlet-to-Neumann
operator given by the harmonic extension over the  cylinder $E \times (0,\infty)$.
 Indeed, let  $  H^1_{L,s}(E\times (0,\infty) ) $ be the set of measurable  functions $w: E \times (0,\infty) \to \R $
with $w\in H^1(E\times (r_1,r_2))$, $0<r_1<r_2<\infty$ and $w=0$ on $ \de E \times (0,\infty) $ such that the following  norm
$$
\|w\|^2_{H^1_{L,s}(E\times(0,\infty)) }=\int_{E \times (0,\infty)  }t^{1-2s}|\n w|^2dxdt<\infty.
$$
In \cite{CT} and \cite{CDDS}, the authors showed that for  any
$g\in \scrH^{-s}(E)$ there exists a unique solution $u\in \scrH^s_0(E)$   to
\be\label{eq:Asvh}
\begin{cases}
 A_su=g\quad \textrm{ in } E,\\
u=0\quad\textrm{ on } \de E.
\end{cases}
\ee
In addition  $u$ is the trace of $w\in H^1_{L,s}(E\times (0,\infty)) $ which is the unique solution to
\be
\begin{cases}
\div(t^{1-2s}\n w)=0\quad \textrm{ in } E\times (0,\infty),\\
w=0 \quad \textrm{ on } \de E \times (0,\infty),\\
-t^{1-2s}\frac{\de w}{\de t}= \k_{N,s}\, g\quad \textrm{ on } E,
\end{cases}
\ee
where $\k_{N,s}$ ($\k_{N,s}=1$ for $s=1/2$) is a constant depending only on $N$ and $s$.
Moreover,
it holds that, with the norm in \eqref{eq:normCT},
\be
\|u\|^2=\k_{N,s}\|w\|^2_{ H^1_{L,s}(E\times (0,\infty) ) }.
\ee

We can compare  the operator $A_s$ with the operator $\calB_s$. For simplicity, we consider the case $s=1/2$.
Assume that $g \in C^\infty_c({E})$ is nonegative and nontrivial and $u$
is a solution to \eqref{eq:Asvh}, which is positive on $E$. Take  $w$ its extension over the cylinder.
 Consider $\calH(\ti{u})$ which is the harmonic extension of $\ti{u}$ in  $\R^{N+1}_+$ given by \eqref{eq:Harmext-intH}.
Clearly
$$
\calH(\ti{u})\geq \ti{w}\quad\textrm{ in } \ov{\R^{N+1}_+}.
$$
It follows from Hopf lemma that
$$
-\frac{\de w}{\de t}> - \frac{\de \calH(\ti{u}) }{\de t}\quad\textrm{ in } E.
$$
Hence
$$
A_{1/2}u>  \calB_{1/2} u\quad\textrm{ in } E .
$$
In particular the operator $A_s$ yields (up to a multiplicative constant) subsolution to the  fractional Laplacian $\Ds$.
This is the reason why the use of  $\calB_s$ is more convenient in this paper.
\end{Remark}
We give here the plan of the paper:
\begin{itemize}
 \item Section \ref{s:NP}: Notations and Preliminaries.
\item Subsection \ref{ss:D-to-N}:  Dirichlet-to-Neumann operator.
\item Section \ref{s:comp-prple}: Comparison and maximum principles.
\item Section \ref{s:nonexist}: Nonexistence of positive  supersolutions.
\item Section \ref{s:exist}: Existence of positive solutions.
\item {Appendix} \ref{a:remainder}, Subsection \ref{ss:Rt}: Remainder term for the fractional Hardy inequality.
\end{itemize}

\section{Notations and Preliminaries}\label{s:NP}
Let $u\in L^2({\R^{N}})$, we will consider
its  Fourier transform
$$
\widehat{u}(\zeta)=\calF(u)(\zeta):=\frac{1}{(2\pi)^{\frac{N}{2}}}\,
\int_{\R^N}e^{-\imath\zeta\cdot{x} }u(x)dx.
$$
For $s>0$, the Sobolev space
$H^s(\R^N)$ is defined as
$$
H^s(\R^N)=\{u\in L^ 2(\R^ N)\,:\, |\zeta|^s\widehat{u}\in   L^2(\R^ N)\}
$$
with  norm
$$
\|u\|_{H^s(\R^N) }= \|\widehat{u}\|_{ L^ 2(\R^ N)}+  \||\zeta|^s \widehat{u}\|_{ L^ 2(\R^ N)}.
$$
We also have by Parseval identity
$$
\|u\|_{H^s(\R^N) }^2=\|u\|_{L^2(\R^N)}^2+ C_{s,N} \int_{\R^N}\int_{\R^N}\frac{|u(x)-u(y)|^2}{|x-y|^{N+2s}}dxdy   .
$$
Let $E$ be a bounded  domain   in $\R^N$ with Lipschitz boundary. For $q>1$, we introduce the space $W^{s,q}(E)$
 defined as the space  of measurable functions $u$  such that the following norm is finite
$$
\|u\|_{W^{s,q}(E) }^q:=\|u\|_{L^q(E)}^q+ \int_E\int_E\frac{|u(x)-u(y)|^q}{|x-y|^{N+qs}}dxdy   .
$$
We define $W^{s,q}_0(E) $ to be the closure of $ C^\infty_c(E)$ with respect to the norm $\|\cdot\|_{W^{s,q}(E) } $.
As a notation convention, we put $H^s(E)=W^{1,2}(E)$ and  $H^s_0(E)=W^{1,2}_0(E)$ which  are Hilbert spaces.\\
It is well known that if $u\in H^1(E)$ then $\ti{u}$,
its null extension outside $E$, is in $H^1(\R^N)$ and $\|u\|_{H^1(E)}=\|u\|_{H^1(\R^N)}$. This is not in general true
for functions in $H^s(E)$ ($s=1/2$ for instance).
We shall define a space of functions in which we recover this defect by imposing integrability of null extensions.\\

The Hardy inequality \eqref{eq:fracHardy}
 suggests the definition of the Hilbert space $\scrD^{s,2}(\R^N)$ which is  the completion of  $C^\infty_c(\R^N)$
with respect to the norm:
\be
v\mapsto \int_{\R^N}|\zeta|^{2s}|\widehat{v}|^2 d\zeta.
\ee
As it will be apparently clear in the remaining of the paper, we  introduce the Hilbert space
\be
\scrH^s_0(E):=\left\{ u\in H^s(E)\,:\,{\ti{u}} \in \scrD^{s,2}(\R^N) \right\},
\ee
where we put here (and hereafter)
$$
\ti{u}=\begin{cases}
        u\quad \textrm{ in } E\\
        0\quad\textrm{ in } \R^N\setminus E.
       \end{cases}
$$
The space $\scrH^s_0(E)$ is endowed with  the norm
\be\label{eq:normscrH}
\| u \|_{\scrH^s_0(E) }^2=\int_{\R^N}|\zeta|^{2s}|\widehat{\ti{u}}|^2 d\zeta=\int_{\R^N}|\zeta|^{2s}|\calF(\ti{u})|^2 d\zeta.
\ee
Note that, since $E$ is bounded, by \eqref{eq:fracHardy}  there exists a constant $C(E)>0$ such that
\be\label{eq:Poincar-scrH}
C(E)\|\ti{u}\|_{ H^s(\R^N) }\leq \| u \|_{\scrH^s_0(E) }\leq \|\ti{u}\|_{ H^s(\R^N) }\quad \forall u \in \scrH^s_0(E).
\ee
Therefore
$$
\scrH^s_0(E)=\{ u\in H^s(E)\,:\, \ti{u} \in  H^s(\R^N)\}.
$$
See for instance [\cite{Gris}, Theorem 1.4.2.2 ], the space $ C^\infty_c(E)$
is  dense in $\scrH^s_0(E)$.\\

\noindent
\textbf{Notations :} For $G$ an open set of $\R^N$, we use the standard notations for weighted
 Lebesgue spaces: $L^p(G; a(x))=\{u:G\to \R\,:\,\int_{G}u^p a(x)dx<\infty \}$ .
 $B^N(0,r)$ is a ball in $\R^N$ centered at 0 with radius $r>0$
and $S^{N-1}= \de B^N(0,1) $.
$\R^{N+1}_+=\{(t,x)\,:\, t>0,\, x\in \R^N\}$. $B^{N+1}_+(0,r)=\R^{N+1}_+\cap B^{N+1}(0,r)$
and $S^{N}_+=\R^{N+1}_+\cap S^{N} $.
If there is no confusion, we will put $ B^N=B^N(0,1)$ and $ B^{N+1}_+=B^{N+1}_+(0,1)$.
The space $W^{1,q}_{0,S}(B^{N+1}_+(0,r); a(x) )=\{u\in W^{1,q}(B^{N+1}_+(0,r); a(x) )\,:\, u=0 \textrm{ on } S^N_+\}$.\\
\noindent
\subsection{ Dirichlet-to-Neumann operator}\label{ss:D-to-N}
 It is well known that the space of Schwartz functions $\calS$ contains
$C^\infty_c({\R^{N}})$ and that $\calF$ is a bijection from $\calS$ into itself.
In  particular ${\Ds\vp}\in C^\infty(\R^N)\cap L^\infty(\R^N) $ for every $\vp\in C^\infty_c({\R^{N}})$.
In fact we have for any $\vp\in C^2_c({\R^{N}})$, see \cite{Sil},
$$
|\Ds\vp(x)|\le C \frac{\|\vp\|_ {C^2_c}}{1+|x|^{N+2s}}\quad \forall x\in \R^N.
$$
This motivates the following:
\begin{Definition}\label{def:fraclap-non-loc}
Let $G$ be an open subset of $\R^N$. Given $u\in \calL^1_s$, the distribution  $\Ds u\in \calD'(G)$
is defined as
$$
\la\Ds u,\vp \ra =\int_{\R^N} u\Ds\vp dx \quad\forall \vp \in C^\infty_c(G).
$$
\end{Definition}
Some recent  results  conserning   $s$-superhamonic functions in the  sense of distributions as above   are in  \cite{Sil}.\\
Consider the Poisson  kernel of $\R^{N+1}_+:=\{(t,x)\,:\, t>0,\, x\in \R^N \}$
\be
P(t,x)=p_{N,s} t^{2s} \frac{1}{(|x|^2+t^2)^{(N+2s)/2}},
\ee
where $p_{N,s}$ is a normalization constant, see \cite{CS} for an explicit value.
Let  $u\in \calL^1_s$, we can define
$$
\bar{u}(t,{x})=P(t,\cdot)*u=p_{N,s}t^{2s}\int_{\R^N}\frac{u(y)}{\left( |y-{x}|^2+ t^2   \right)^{\frac{N+2s}{2}}   }
dy\quad\forall (t,x)\in \R^{N+1}_+.
$$
It turns out that
$$
\div(t^{1-2s}\n\bar{u})=0\quad \R^{N+1}_+.
$$
Therefore $\bar{u}$ is smooth in $\R^{N+1}_+$.
Moreover if
 $u$ is regular in a neighborhood of some point $x_0$ then
$$
\lim_{t\to 0} \bar{u}(t,{x_0})\to u(x_0) .
$$
By an argument of \cite{CSilv}, we have that
\be
-\lim_{t\to 0} t^{1-2s}\frac{\de \bar{u}}{\de t}(t,x_0)= \k_{s} \Ds u(x_0),
\ee
where the constant $\k_{s}$ is explicitly computed in  \cite{CS}:
\be\label{eq:kappas}
k_s=\frac{\G(1-s)}{2^{2s-1}\G(s)}.
\ee
%
For any $w\in H^s(\R^N)$, we denote by $\calH(w) $ its unique harmonic extension over $\R^{N+1}_+$.
Namely (see for instance \cite{CSilv}, \cite{CS}) $\calH(w)\in H^1(\R^{N+1}_+;t^{1-2s})$ and
\be\label{eq:Harmext}
\begin{cases}
\div(t^{1-2s}\n \calH(w))=0\quad \textrm{ in }\R^{N+1}_+,\\
\calH(w) = w \quad \textrm{ on }\R^N,\\
-t^{1-2s}\frac{\de \calH(w)}{\de t}= \k_{s} \Ds w \quad   \textrm{ on } \R^N.
\end{cases}
\ee
In particular if $w\in C^2_c(\R^N)$ then $\calH(w)=P(t,\cdot)*w $.
In addition  one can  check (see \cite{CSilv} ), using integration by parts  and the Parseval identity, that
\be\label{eq:ncaHeqns2}
\int_{\R^{N+1}_+}t^{1-2s}|\n \calH(w)|^2dxdt=\k_{s}\int_{\R^N}|\zeta|^{2s}\widehat{w}d\zeta
= \k_{s}\int_{\R^N}|(-\D)^{s/2}w|^2 dx.
\ee
Therefore from the  definition of the space $\scrH^s_0(E) $, we have
\be\label{eq:scrHeqH1calH}
\k_{s}\|v\|_{\scrH^s_0(E)}^2= \int_{\R^{N+1}_+}t^{1-2s}|\n \calH(\ti{v})|^2dxdt\quad\forall v\in \scrH^s_0(E),
\ee
where as usual $\ti{v} $ is the null extension of $v$ outside $E$.
\\

We now  introduce a Dirichlet-to-Neumann  operator $\calB_s$ defined on
  $\scrH^s_0(E)$.
\begin{Proposition}\label{prop:DtN}
Let $E$ be a bounded open set with Lipschitz boundary. Denote by $\scrH^{-s}(E)$  the dual of $\scrH^s_0(E)$.
Then the mapping
   $\calB_s: \scrH^s_0(E) \to\scrH^{-s}(E) $ given by
$$
 \la \calB_s v,\vp\ra_{\scrH^{-s}(E), \scrH^s_0(E) } =-\k_s^{-1} \int_{\R^N}\lim_{t\to 0} t^{1-2s}
\frac{\de \calH(\ti{v})}{\de t}\ti{\vp} dx\quad \forall v, \vp\in \scrH^s_0(E)
$$
is a linear  isometry. In addition
 for any $v\in \scrH^s_0(E)$ we have
\be\label{eq:D-to-N}
 \calB_s v= \Ds\ti{v}\quad \textrm{ in }  \calD'(E) .
\ee
\end{Proposition}
\proof
By definition for any  $v\in  \scrH^s_0(E)$, $\ti{v} \in H^s(\R^N)$ thus the operator $\calB_s$ is well defined and linear.
 Consider $\calH(\ti{v})$ which satisfies \eqref{eq:Harmext}.
Then integration by parts yields for every $\vp\in\scrH^s_0(E)$
\begin{eqnarray*}
 \int_{\R^{N+1}_+}t^{1-2s}\n \calH(\ti{v})\cdot\n \calH({\ti{\vp}})dxdt
&=&\int_{\R^N} \lim_{t\to 0} t^{1-2s}\frac{\de \calH(\ti{v})}{\de t}\ti{\vp} dx.
\end{eqnarray*}
This, \eqref{eq:scrHeqH1calH}  and H\"{o}lder inequality imply that
$$
\la \calB_s v,\vp \ra_{\scrH^{-s}(E), \scrH^s_0(E) }\leq \|v\|_{\scrH^s_0(E)}^2\|\vp\|_{\scrH^s_0(E)}^2
$$
while
$$
\la \calB_s v,v  \ra_{\scrH^{-s}(E), \scrH^s_0(E) }=\|v\|_{\scrH^s_0(E)}^2.
$$
Hence
$$
\|\calB_s v\|_{\scrH^{-s}(E) }=\|v\|_{\scrH^s_0(E)}.
$$
On the other hand, for any $\vp\in C^\infty_c(E)$, we have by   integration by parts
\begin{eqnarray*}
 \int_{\R^{N+1}_+}t^{1-2s}\n \calH(\ti{v})\cdot\n \calH({\ti{\vp}})dxdt
&=&\int_{\R^N} \lim_{t\to 0} t^{1-2s}\frac{\de\calH(\ti{v})}{\de t}\vp dx\\
&=&\int_{\R^N} \lim_{t\to 0} t^{1-2s}\frac{\de\calH(\vp)}{\de t}\ti{v} dx\\
&=&\k_s\int_{\R^N}\ti{v} \Ds {\vp}  dx.
\end{eqnarray*}
This means that
$$
\calB_s v= \Ds\ti{v} \quad\textrm{ in } \calD'(E).
$$
\QED
We turn to the characterization of the space $\scrH^s_0(E)$.
As suggested with the fact that $\calH(\ti{v})\in H^1(\R^{N+1}_+;t^{1-2s}) $ for every $ v\in \scrH^s_0(E)$,
we have the converse:
\begin{Proposition}\label{prop:H1zT-eqscrH}
Let $E$ be a bounded open set with Lipschitz boundary.
Define
$$
H^1_{0,T}(E;t^{1-2s})=\left\{w\in H^1(\R^{N+1}_+;t^{1-2s})\,:\,   w\Big|_{\R^N}\equiv0 \textrm{ on } \R^N\setminus E \right \}
$$
and
$$
H^s_{0,T}(E)=\left\{u\in \calD^{s,2}(\R^{N})\,:\,   u\equiv0 \textrm{ in } \R^N\setminus E \right \}.
$$
 We have the following equalities:
\be\label{eq:cahrac-scrHsz}
\scrH^s_0(E)= \left\{u\Big|_{E}\,:\, u\in H^s_{0,T}(E)\right\}=\left\{w\Big|_{E}\,:\, w\in H^1_{0,T}(E;t^{1-2s})\right\}.
\ee
In particular
$$
\Ds u=\calB_s\check{u} \quad\textrm{ in } \calD'(E),\quad \forall u\in H^s_{0,T}(E),
$$
where $ \check{u}=u \Big|_{E} $.
\end{Proposition}
\proof
The first equality in \eqref{eq:cahrac-scrHsz} is immediate  by definition.
The second equality is a consequence of the trace embedding theorem.
Indeed, take $w\in H^1_{0,T}(E;t^{1-2s})$. Then the null extension of
${w\Big|_{E} }$ outside $E$ is nothing but $w$ which belongs to $ H^s(\R^N)$  and in addition
$
 \|w\|_{H^{s}(E)}\leq  \|w\|_{H^{s}(\R^N)} .
$
\QED
Summarizing, we state the following
\begin{Proposition}\label{prop:Di-to-Neu}
Pick $g\in \scrH^{-s}(E)$. Let $v\in\scrH^s_0(E)$ (given by the Lax-Miligram theorem) be the unique  solution to
$$
\calB_s v=g\quad \textrm{ in } E.
$$
Let $w\in H^1(\R^{N+1}_+;t^{1-2s})$ solve the mixed problem
$$
\begin{cases}
\div(t^{1-2s}\n w)=0\quad \textrm{ in }\R^{N+1}_+,\\
w = 0 \quad \textrm{ on }\R^N\setminus E,\\
-t^{1-2s}\frac{\de w}{\de t}= \k_{s} g \quad \textrm{ on } E.
\end{cases}
$$
Then $v=w$ in $E$; for any $\vp\in \scrH^s_0(E)$
\begin{eqnarray*}
\la \calB_s v,\vp \ra_{\scrH^{-s}(E), \scrH^s(E) }
&= & \int_{\R^N}|\zeta|^{2s}\calF(\ti{v})\calF(\ti{\vp})\\
&=&\int_{\R^N}(-\D)^{s/2} \ti{v}(-\D)^{s/2}\ti{\vp} dx \\
&=&\la g,\vp \ra_{\scrH^{-s}(E), \scrH^s(E) }\\
&=&\k_s^{-1}\int_{\R^{N+1}_+}t^{1-2s}\n w\cdot\n \calH(\ti{\vp})dxdt\\
&=&\k_s^{-1}\int_{\R^{N+1}_+}t^{1-2s}\n \calH(\ti{v}) \cdot\n \calH(\ti{\vp})dxdt
\end{eqnarray*}
 and thus
$$
\k_s\|v\|^2_{\scrH^s_0(E)}=\int_{\R^{N+1}_+} t^{1-2s}|\n w|^2dxdt.
$$
\end{Proposition}

We can extend the above in unbounded domains:
\begin{Remark}
Here we consider $E$ any open subset of $\R^N$ with Lipschitz boundary. Define
\be
\scrH^s(E):=\left\{ u\in H^s(E)\,:\,{\ti{u}} \in H^1(\R^N) \right\},
\ee
where as usual $\ti{u}$ stands for the null extension of $u$ outside $E$.
We have that $ C^\infty_c(E)$ is dense in $\scrH^s(E)$, see \cite{Gris}.\\
By similar arguments, we have that the operator
$$
\bar{\calB_s}(v)= -\k_s^{-1}t^{1-2s}\frac{\de  \calH(\ti{v})}{\de t}+ v
$$
is a linear isometry form $\scrH^s(E)\to (\scrH^s(E))' $, where $(\scrH^s(E))' $ is the dual of $\scrH^s(E)$.
\end{Remark}
%

\section{Comparison and maximum principles}\label{s:comp-prple}
Unless otherwise stated, $E$ is a bounded Lipschitz open set  of $\R^N$.
We have the following technical result which will be useful in the sequel.
\begin{Lemma}\label{lem:E_n}
Let $E_n$ be  a sequence of Lipschitz open sets  such that $E_n\subset\subset E_{n+1}$ and $\cup_{n=1}^{\infty} E_n=E$.
Let $g_n\in L^2(E)$ such that $g_n\to g$ in $L^2(E)$. Consider $v_n\in\scrH^s_0(E_n)$ solution to
$$
\calB_s v_n=g_n\quad \textrm{ in } E_n.
$$
If $v\in\scrH^s_0(E)$ is the unique solution to
$$
\calB_s v=g\quad \textrm{ in } E
$$
then $\ti{v_n}\to v$ in $L^2(E)$.
\end{Lemma}
\proof
Observe that $\calH(\ti{v_n})\in H^1_{0,T}(E; t^{1-2s})$ thus by Proposition \ref{prop:H1zT-eqscrH} $\ti{v_n}\in  \scrH^s_0(E)$. In addition we have by
Hardy and H\"{o}lder inequality
$$
\|\ti{v_n}\|_{\scrH^s_0(E)  }=\|v_n\|_{\scrH^s_0(E)  }\leq C(E) \|g_n\|_{L^2(E)}.
$$
Therefore $\ti{v_n}$ is bounded. By assumption it   converges weakly to $v$ in $\scrH^s_0(E)$ and strongly in $L^2(E)$
because $C^\infty_c(E)$ is dense in $ \scrH^s_0(E)$.
 \QED


The following maximum principle can be found in [\cite{CDDS} Lemma 2.4] or in  \cite{FKS}.
\begin{Lemma}\label{lem:smmp}
Let  $E$ be a bounded Lipschitz domain of $\R^N$.
Let  $v\in \scrH^s_0(E)$, $v\geq 0$ such that
$$
\calB_s v \geq0 \quad \textrm{ in } E.
$$
If $v  \neq0$  then  for any compact set $K\subset E$
$$
\textrm{ess}\inf_K {v}>0.
$$
\end{Lemma}


\begin{Lemma}\label{lem:mmp1}
Let $g\in L^2(E)$, $g\geq0$ and let $w\in L^{1}_{loc}(\ov{\R^{N+1}_+ })$,    such that
$$
\int_{\R^{N+1}_+}t^{1-2s}|\n w|^2dxdt<\infty
$$
and
\be\label{eq:weq-mp}
\int_{\R^{N+1}_+}t^{1-2s}\n w\cdot\n\phi dxdt+c\int_{E}w\vp dxdt\geq \k_s\int_{E}g\phi dx
\ee
for every  nonegative $\phi\in H^1_{0,T}( E; t^{1-2s})$, where $c\in\R_+$. Assume that $w\geq0$ on $\R^{N}\setminus E$.
Then $w\geq   0$ in $\ov{\R^{N+1}_+}$.
\end{Lemma}
\proof
Test  \eqref{eq:weq-mp} with $\max(-w,0)\in H^1_{0,T}( E; t^{1-2s})$.
\QED

\begin{Lemma}\label{lem:min-solH10-non-loc}
Let $c\in\R_+$ and let    $u\in
\calL^1_s$, $u\geq 0$ and $g\in L^2(E)$  such that
\be\label{eq:Dsugqgmi}
\Ds u+cu\geq g \quad \textrm{ in } \calD'(E).
\ee
 Let $v\in \scrH^s_0(E )$ solves
\be
\calB_s v+cv=g \quad \textrm{ in } E.
\ee
Then
 $$
u\geq  v\quad\textit{ in $E$.}
$$
\end{Lemma}
\proof
Recall that \eqref{eq:Dsugqgmi} is equivalent to
\be\label{eq:defvwksp}
\int_{\R^N}u\Ds\vp dx \geq \int_{E }g\vp dx-c \int_{E }u\vp dx \quad\forall \vp\in C^\infty_c(E), \,\vp\geq0.
\ee
Denote by $\rho_n$  the
standard mollifier (which is symmetric: $\rho_n(-x)= \rho_n(x)$) and put $u_n=\rho_n*{u}$.\\
\textbf{Claim: }
 for any $\vp\in C^\infty_c(\R^N) $
\be\label{eq:convol}
\int_{\R^N}\Ds u_n\vp=\int_{\R^N}u\Ds(\rho_n*\vp).
\ee
It  is easy to check  using Fubini's theorem and the symmetry of $\rho_n$ that
\be\label{eq;integ-part}
 \int_{\R^N}\Ds u_n\vp dx=\int_{\R^N} u_n {\Ds\vp}dx=\int_{\R^N}    {u} \rho_n*{\Ds\vp}dx.
\ee
Now  we notice that, in $\R^N$,
$$
 \rho_n*{\Ds\vp}=\calF(\calF(\rho_n*{\Ds\vp}))=\calF( |\zeta|^{2s}(\calF(\rho_n)\calF(\vp)))= \Ds(\rho_n*\vp).
$$
Using this in \eqref{eq;integ-part},  we get \eqref{eq:convol} as claimed.\\
 Let $E_n:=\{x\in E\, :\, \textrm{dist}(x,\de E)>1/n\}$.
We deduce from \eqref{eq:defvwksp} and \eqref{eq:convol} that for all $\vp\in C^\infty_c(E_n)$ and $\vp\geq0$
\begin{eqnarray*}
\int_{\R^N }\Ds u_n\vp dx=\int_{\R^N }u\Ds(\rho_n*\vp) dx &\geq &\int_{E }g(\rho_n*\vp)dx-c\int_{E}u (\rho_n*\vp)dx \\
&=& \int_{E }(\rho_n*g)\vp dx-c\int_{E }(\rho_n*u)\vp dx .
\end{eqnarray*}
We conclude that
\be
\Ds u_n(x)+c u_n(x)\geq \rho_n*g(x) =:g_n(x)\quad\textrm{ for every   } x\in E_n.
\ee
We let  ${w}_n(t,x)=P(t,\cdot)*u_n(x)$ be the harmonic extension of $u_n$ via the Poisson kernel so that 
\be
\begin{cases}
\div(t^{1-2s}\n w_n)=0\quad\R^{N+1}_+ ,\\
w_n= u_n \quad  {\R^N}.
\end{cases}
\ee
It turns out that
\be
-t^{1-2s}\frac{\de w_n}{\de t}+c w_n=\k_s \Ds u_n+c u_n\geq \k_s g_n\quad \textrm{ on  }E_n 
\ee
and in addition $t^{1-2s}|\n w_n|^2\in L^1_{loc}(\ov{\R^{N+1}_+})$.
Let  $v_n\in \scrH^s_0(E_n) $ be the solution to $$\calB_s v_n+c v_n=g_n \quad\textrm{ in } E_n.$$
We take a large  $R>0$ so that $B^N(0,R)$ contains $E$ and  we let  $v_{n,R}\in W^{1,2}_{0,S}(B^{N+1}_+(0,R);t^{1-2s} )$ 
be the unique solution (obtained by minimization) to the problem
\be
\begin{cases}
\div(t^{1-2s}\n v_{n,R})=0 \quad B^{N+1}_+(0,R) ,\\
-t^{1-2s}\frac{\de v_{n,R}}{\de t}+c  v_{n,R}= \k_sg_n \quad  E_n,\\
v_{n,R}=0\quad B^{N}(0,R)\setminus E_n.
\end{cases}
\ee
By extending $v_{n,R}$ to be zero outside $ \ov{B^{N+1}_+(0,R)}$, it is standard to show that $  v_{n,R}\to v_n$ as $R\to\infty$ in $ H^1(\R^{N+1}_+;t^{1-2s})$. 
Since  ${w}_n\geq  v_{n,R}$ by  Lemma \ref{lem:mmp1}, it follows that, sending $R\to \infty$,  $ {w}_n\geq  v_n$  in $\R^N$.
 In particular $u_n\geq v_n$ in $E_n$. 
By Lemma \ref{lem:E_n}, $ \ti{v_n} \to v $ in $L^2(E)$
and the proof is complete.
 \QED
%

We recall the  definition of the  $s$-capacity of a compact set $A\subset E$:
\be
C_s(A)=\inf_{\phi \in C^\infty_c(E) }\{ \|\phi \|^2_{\scrH^{s}_0(E) }\,:\, \vp \geq 1 \textrm{ in a neighborhood of $A$}\}.
\ee
Note that if $ C_s(A)=0$ then $ |A|=0$ by Poincar\'e inequality (see \eqref{eq:Poincar-scrH}).
 We have the following comparison result modulo small sets.
\begin{Lemma}\label{lem:min-solH10-loc}
 Let $A$ be a compact subset of $E$ with $C_s(A)=0$. Let $u\in\calL^1_s$, $c\in\R_+$ and $g\in L^2(E)$ such that
 $$
\Ds u +c u \geq g \quad \textrm{ in } \calD'(E\setminus A).
$$
Let $v\in\scrH^s_0(E)$ solve
$$
\calB_s v +cv=g\quad \textrm{ in } E.
$$
Then $u\geq v$ in $E$.
\end{Lemma}
\proof
Let  $A_\e$ be  a smooth open $\e$-neighborhood of $A$ compactly contained in $E$.
 Define $D_\e=\{x\in D\,:\,\textrm{dist}(x,\de (D\setminus A_\e))>\e\}$.
 It is clear that
$$
\Ds u +cu \geq g\quad \textrm{ in } \calD'({D}_\e).
$$
 Consider $v_\e\in\scrH^s_0({D}_\e)$
solving $$\calB_s v_\e +cv_\e=g \quad \textrm{ in }{D}_\e.$$
 By Lemma \ref{lem:min-solH10-non-loc} we have $u\geq v_\e$ in $D_\e$.
The same argument as in the proof of Lemma \ref{lem:E_n}  yields $\ti{v_\e}\in \scrH^{s}_0(E) $ for every $s\in(0,1)$
and it is bounded. Hence it converges weakly to some function $w$ in  $\scrH^{s}_0(E)$ and strongly in $L^2(E)$.
In particular $u\geq w$. Moreover for any $\vp\in C^\infty_c(E\setminus A)$, we can choose $\e>0$ so
 small that $supp\vp$ is contained in $D_\e $ thus
 taking the limit as $\e\to0$, we get
$$
\la w,\vp\ra_{\scrH^{s}_0(E) } +c\int_{E}w\vp=\int_{E}g\vp dx\quad\forall \vp\in C^\infty_c(E\setminus A).
$$
From this equality, to conclude the proof (that is $v=w$),
it suffices to show that $ C^\infty_c(E\setminus A)$ is dense in $C^\infty_c(E) $
with the $\scrH^{s}_0(E) $-norm because $w\in \scrH^{s}_0(E)$.\\
Since $C_s(A)=0$, there exists a sequence $\psi_n \in C^\infty_c(E)$ such that $\psi_n \geq 1$ in a neighborhood of $A$ and
in addition
\be\label{eq:psi_nto0}
 \|\chi_n\|^2_{\scrH^{s}_0(E) }\leq \|\psi_n\|_{\scrH^{s}_0(E) }\to0,
\ee
where $\chi_n=\min(\psi_n,1) $.
Now take any $\phi\in C^\infty_c(E)$  and note that  $(1-\chi_n)\phi\in C^\infty_c(E\setminus A)$ and
 moreover $(1-\chi_n)\phi\to \phi$ in $\scrH^{s}_0(E) $
 by  \eqref{eq:psi_nto0}. This concludes the proof.
\QED
%
We shall define a new space which is more convenient when  dealing with the Hardy potential.
Namely, we assume that there exists $b\in L^1_{loc}(E)$  and a constant $C>0$ such that
\be\label{eq:coerciv}
\|\vp\|^2_{ \scrH^s_{0} (E)}- \int_{E}b(x)\vp^2dx\geq C\int_{E}\vp^2dx \quad\forall \vp\in C^{\infty}_c( E).
\ee
\begin{Definition}\label{deqf:scrH0b}
Let $b\in L^1_{loc}(E)$ so that \eqref{eq:coerciv} holds .
The  Hilbert space $\scrH^s_{0,b} (E) $ is  the completion of $C^\infty_c(E)$ with respect to the scalar product
$$
\la\vp,\phi\ra_{\scrH^s_0 (E) }-  \int_{E}b(x)\vp\phi dx\quad\forall \vp,\phi\in C^{\infty}_c(E).
$$
\end{Definition}
Note that the Lax-Miligram theorem implies that for any $f\in L^2 (E) $, there exits a unique
solution to the problem
\be\label{eq:PBgb}
 \begin{cases}
 \calB_s v-b(x)v=f\quad \textrm{ in } E,\\
   v\in \scrH^s_{0,b}(E)  ,
\end{cases}
 \ee
in the sense that for all $\phi \in \scrH^s_{0,b}(E)$
$$
\la v,\phi\ra_{\scrH^s_0 (E) }-  \int_{E}b(x)v\phi dx=\int_{E}f\phi dx.
$$
\begin{Remark}\label{rem:scrHdepscrHz}
Let $\e>0$. Put $d_\e(x)=b(x)(1-\e)$. Then  $ \scrH^s_{0}(E)=  \scrH^s_{0,d_\e}(E) $
by Propositon \ref{prop:H1zT-eqscrH}.
This holds true because if $v\in \scrH^s_{0,d_\e}(E)$ then by  \eqref{eq:coerciv} we have
$\ti{v}\in H^s(\R^N)$. By similar argument $\scrH^s_{0}(E)=  \scrH^s_{0,b}(E)$ if $b\in L^\infty(E)$.
\end{Remark}
\begin{Lemma}\label{lem:Complem}
Let $A$ be a compact subset of $E$ with $C_s(A)=0$.
Let $b\in L^1_{loc}(E )$ such that \eqref{eq:coerciv} holds.
Suppose that $u\in \calL^1_s$ with  $u,b\geq 0$ and $f\in  L^2(E)$, $f\geq 0$ such that
\be\label{eq:Dsugqg}
\Ds u- b(x) u \geq f\quad\textrm{ in } \calD'(E\setminus A).
\ee
Let $v\in\scrH^s_{0,b}(E) $ be the unique solution to
$$
\calB_s v-b(x)v=f\quad \textrm{ in } E.
$$
Then
$$
u\geq  v \quad \textrm{ in }E.
$$
\end{Lemma}
\proof
\textbf{Step 1:}\textit{ We first prove the result  if  $b\in L^\infty(E)$.}\\
 We  let  $v_0\in \scrH^s_0(E)$ solving
$$
 \calB_s v_0= f\quad\textrm{ in $E$}.
$$
Then $ 0\leq v_0\leq u$ in $E$ by
Lemma \ref{lem:min-solH10-loc} and because $f\geq 0$.
We define inductively the sequence $v_n\in \scrH^s_0(E)$ by
$$
 \calB_s v_1= {b}(x) v_{0}+f\quad\textrm{in $E$},\quad\quad
 \calB_s v_n= {b}(x) v_{n-1}+f \quad\textrm{in $E$}.
$$
Since ${b}\geq 0$, we have  $\Ds u\geq  {b}(x)v_0+ f$ in $\calD'(E\setminus A) $. Thus   using  once again Lemma \ref{lem:min-solH10-loc}, we  obtain
$ v_0\leq v_{1}\leq u $ in $E$. By induction, we  have  
$$
v_0\leq v_1\leq \dots\leq  v_{n}\leq u\quad\textrm{  in $E$}\quad \forall n\in\N.
$$
Since $v_{n-1}\leq v_n$ in $E$, we have
$$
\|v_n\|^2_{\scrH^s_{0}(E)}-\int_{E}b(x)|v_n|^2\leq \int_{E}f(x)v_ndx.
$$
By H\"{older} inequality and \eqref{eq:coerciv} (see Remark \ref{rem:scrHdepscrHz}) $v_n$ is bounded in $\scrH^s_0(E)$.
We conclude  that $v_n \rightharpoonup v$  in $\scrH^s_0(E)$ as $n\to \infty$ which is the unique solution to 
$$
\calB_s v= {b}(x) v+f\quad\textrm{in $E$}.
$$
Since $v_n\to v$ in $L^2(E)$, we get  $v\leq u$ in $E$.

\bigskip 
\noindent
\textbf{Step 2:} \textit{Conclusion of the proof.}\\
We put $b_k(x)=\min(b(x),k)$ for every $k\in\N$. We consider 
 ${v}^k\in \scrH^s_0(E)$  be the unique solution to 
\begin{equation}\label{eq:tvkstf}
 \la {v}^k,\f\ra_{\scrH^s_{0}(E)}-\int_{E}\min\left\{b(x),{k}\right\} {v}^k\f=\int_{E}f\f\quad\forall\f\in C^\infty_c(E).
\end{equation}
Thanks to   \textbf{Step 1}, we have $v^k\leq u$ in $E$.\\
Next, we check that  such a sequence  ${v}^k$, satisfying \eqref{eq:tvkstf}, 
  converges to $v$ in $L^2(E)$ when $k\to\infty$. Indeed, we have 

\begin{eqnarray*}
\| {v}^k\|^2_{{\scrH^s_{0,b}(E)} }
&\leq& \| {v}^k\|^2_{{\scrH^s_0(E)} }-\int_{E}\min\{b(x),k\}~\!|{v}^k|^2~dx\\
&=&
\int_{E} f{v}^k~dx\le C\|{v}^k\|_{\scrH^s_{0,b}(E)}
\end{eqnarray*}
by H\"older inequality and by (\ref{eq:coerciv}), where the constant $C$
depends on $f$ and $E$ but not on $k$.
Therefore the sequence ${v}^k$ is bounded in ${\scrH^s_{0,b}(E)}$. We conclude that  there exists
$\ti{v}\in {\scrH^s_{0,b}(E)}$ such that, for a subsequence, ${v}^k \rightharpoonup \ti{v}$ in
${\scrH^s_{0,b}(E)}$.
Now by \eqref{eq:tvkstf},  we have
$$
\la {v}^k,\f\ra_{{\scrH^s_{0,b}(E)}}+\int_{E}\left(b(x)-
\min\{b(x),k\}\right){v}^k\vp=\int_{E}f\vp.
$$
Since for every $k\geq 1$ and any $\vp\in C^\infty_c(E) $
 $$\left|\left(b(x)- \min\{b(x),k\}\right){v}^k\vp\right|\leq\left(b(x)- \min\{b(x),k\}\right)u|\vp|
\leq 2b(x)u|\vp|\in L^1(E),$$
 the dominated convergence theorem implies
that
\begin{equation}\label{eq:vstw}
\la \ti{v},\f\ra_{\scrH^s_{0,b}(E)}=\int_{E} f \vp\quad\textrm{for any $
\vp\in C^\infty_c(E)$.}
\end{equation}
We therefore have  that $\ti{v}=v$ by uniqueness.
By \eqref{eq:vstw}, we have
\begin{eqnarray*}
 \| {v}-{v}^k\|^2_{{\scrH^s_{0,b}(E)}}&=&
\| {v}^k\|_{{\scrH^s_{0,b}(E)}}^2- \la {v},{v}^k\ra_{{\scrH^s_{0,b}(E)}}
+\la {v},{v}-{v}^k\ra_{{\scrH^s_{0,b}(E)}}\\
&=&\| {v}^k\|_{{\scrH^s_{0,b}(E)}}^2-\int_{E} f {v}^k+\la {v},{v}-{v}^k\ra_{{\scrH^s_{0,b}(E)}}\\
&\leq&  \| {v}^k\|^2_{\scrH^s_0(E)}-\int_{E}\min\{b(x),k\}~\!|{v}^k|^2~dx -\int_{E} f {v}^k+\la {v},{v}-{v}^k\ra_{{\scrH^s_{0,b}(E)}}\\
&=&   \la {v},{v}-{v}^k\ra_{{\scrH^s_{0,b}(E)}} .
\end{eqnarray*}
We thus obtain 
$$
C(E)\,\int_{\O}|{v}-{v}^k|^2~dx\leq \la {v},{v}-{v}^k\ra_{{\scrH^s_{0,b}(E)}}\to 0
$$
by \eqref{eq:coerciv}.
Hence  ${v}^k\to {v}$ pointwise and
thus ${v}\leq u$ in $\O$.
\QED
\begin{Remark}\label{rem:complem-c}
 The same result as in Lemma \ref{lem:Complem} holds if we assumed
the coercivity that there exist  constants $C,c>0$ such that
for all $\vp\in C^{\infty}_c( E)$
\be
\|\vp\|^2_{\scrH^s_0 (E) }+c \int_{E}\vp^2 dx-
\int_{E}b(x)\vp^2 dx\geq C\int_{E}\vp^2 dx .
\ee
%
\end{Remark}
We close this section with the following useful lemma and its immediate consequence. Its counterpart, for $s=1$,
is in \cite{FaMu-ne}.
\begin{Lemma}\label{lem:AP}
Let  $E$ be a  bounded Lipschitz domain of $\R^N$.
Let $A$ be a compact subset of $E$ with $C_s(A)=0$.
 Let $u\in\calL^1_s$,  $b\in L^1_{loc}(E)$ and $u,b>0$. Assume that
\be\label{eq:Dsugqg}
\Ds u\geq b(x) u\quad\textrm{ in }\calD'(E\setminus A).
\ee
Then
\be\label{eq:AP}
\int_{\R^{N+1}_+ }t^{1-2s}|\n \calH(\vp) |^2dxdt=\k_s\|\vp\|^2_{\scrH^s_0(E)}\geq
\k_s \int_{E}b(x)\vp^2dx\quad\forall \vp\in C^{\infty}_c(E).
\ee

\end{Lemma}
\proof
Put $g_k(x):=\min( b(x) u ,k )>0$  for integers $k\geq1$. Let $v_k\in \scrH^s_0(E)$ be the solution to
$$
\calB_s v_k=g_k\quad\textrm{ in } E.
$$
By Lemma \ref{lem:smmp}, we have $\frac{1}{v_k}\in L^\infty_{loc}(E ) $
and by the standard maximum principle $\calH(\ti{v_k})>0 $.
Moreover by Lemma \ref{lem:min-solH10-loc}, we have
\be\label{eq:ugeqvk}
u\geq v_k>0 \quad \textrm{ in } E.
\ee
 Let $\vp\in C^{\infty}_c( E)$. Put $V_k= \calH(\ti{v_k})$ and $ V_k^\e= V_k+\e$,  for $\e>0$.
Set $\psi=\frac{\calH(\vp)}{V_k^\e}$ so that
 $V_k^\e\psi^2\in H^1_{0,T}(E; t^{1-2s})$.
Simple computations show that
$$
|\n \calH(\vp)|^2=|V_k^\e\n\psi|^2+\n V_k^\e\cdot \n({V_k^\e}\psi^2)= |V_k^\e\n\psi|^2+\n V_k\cdot \n({V_k^\e}\psi^2).
$$
Thus using \ integration by parts we have
\begin{eqnarray*}
 \int_{\R^{N+1}_+ }t^{1-2s}|\n \calH(\vp) |^2dxdt&\geq&  \int_{\R^{N+1}_+ }t^{1-2s} \n V_k\cdot \n(V_k^\e\psi^2)dxdt\\
&= &\int_{E}{g_k}\frac{\vp^2}{(v_k+\e)^2}dx.
\end{eqnarray*}
Take the limit as $\e\to0$ to get
$$
 \int_{\R^{N+1}_+ }t^{1-2s}|\n \calH(\vp) |^2dxdt\geq \k_s \int_{E}\frac{g_k}{v_k}{\vp^2}dx
$$
by Fatou's lemma.
By \eqref{eq:ugeqvk}, we infer that
$$
 \int_{\R^{N+1}_+ }t^{1-2s}|\n \calH(\vp) |^2dxdt\geq \k_s \int_{E}\frac{g_k}{u}{\vp^2}dx.
$$
Again  by  Fatou's lemma, inequality \eqref{eq:AP} follows immediately by taking $k\to+\infty$.
\QED
The following result appeared in \cite{BG} in the case $s=1$.
\begin{Theorem}\label{th:necess}
 Let $E$ be a bounded Lipschitz domain of $\R^N$ with $0\in E $, $N>2s$. Then there is no nonnegative and nontrivial $u\in \calL^1_s$ satisfying
$$
\Ds u\geq \g |x|^{-2s} u\quad \textrm{ in } \calD'(E\setminus\{0\}),
$$
with $\g> \g_0$.
\end{Theorem}
\proof
Note  that $C_s(\{0\})=0$ provided $N>2s$ (see \cite[p. 397]{Maz}).
If such $u$ exits then $u>0$ in $E$ by the maximum principle thus Lemma \ref{lem:AP}
 contradicts the sharpness of the Hardy constant $\g_0$.
\QED


\section{Nonexistence of positive  supersolutions}\label{s:nonexist}
We start with the following

\begin{Lemma}\label{lem:grstate}
 For every $\a\in(-\frac{N}{2}-s,\frac{N}{2}-s)$, put $\vartheta_\a({x})= |{x}|^{\frac{2s-N}{2}+\a}$. Then
$$
\Ds \vartheta_\a=\g_\a |{x}|^{-2s}\,\vartheta_\a\quad \textrm{ in } { \R^N\setminus\{0\}},
$$
where
\be\label{eq:defgamal}
\g_\a=2^{2s} \frac{\Gamma\left(\frac{N+2s+2\a}{4}\right)}{\Gamma\left(\frac{N-2s-2\a}{4}\right)} 
\frac{\Gamma\left(\frac{N+2s-2\a}{4}\right)}{\Gamma\left(\frac{N-2s+2\a}{4}\right)}
.
\ee
For $\a\geq0$, the  function $\a\mapsto \g_\a$ is continuous and  decreasing.\\
 There exists a positive function $\Upsilon_\a\in C^\b\left(\ov{\R^{N+1}_+}\setminus\{0\}\right)$ such that
\be\label{eq:v-al_def}
\begin{cases}
\div(t^{1-2s}\n \Upsilon_\a)=0\quad \textrm{ in } \R^{N+1}_+\\
\Upsilon_\a=\vartheta_\a\quad  \textrm{ on }\de \R^{N+1}_+\setminus\{0\}\\
-t^{1-2s}\frac{\de \Upsilon_\a }{\de t}=\k_s  \Ds \vartheta_\a
=\k_s\g_\a |{x}|^{-2s}\,\vartheta_\a\quad\textrm{ on }\de \R^{N+1}_+\setminus\{0\}.
\end{cases}
\ee
Moreover if $\a>0$ then $|\n \Upsilon_\a|\in L^2(B^{N+1}_+(0,R);t^{1-2s})$ for every $R>0$.
\end{Lemma}
\proof
Note that $ \vartheta_\a\in \calL^1_s$.
The Fourier transform of radial functions (see [\cite{SW} Theorem 4.1]) yields
$$
\calF(\vartheta_\a)(\rho)=\rho^{\frac{1-N}{2}}\int_0^\infty(r\rho)^{\frac{1}{2}} J_{\frac{N-2}{2} }(r\rho)\vartheta_\a r^{\frac{N-1}{2} }dr,
$$
where $  J_{\frac{N-2}{2} }$ is the Bessel function. Then we have
\begin{eqnarray*}
\calF(\vartheta_\a)(\rho)&=&\rho^{\frac{-N}{2}-s-\a} \int_0^\infty(r\rho)^{s+\a} J_{\frac{N-2}{2} }(r\rho)d(r\rho)\\
&=&{m_\a} \rho^{\frac{-N}{2}-s-\a},
\end{eqnarray*}
where  $$ m_\a=2^{s+\a}\frac{\G\left(\frac{N+2s+2\a}{4  }\right)}{\G\left(\frac{N-2s-2\a}{4  }\right)}.$$
Now we notice that ($\g_\a=m_\a m_{-\a}$)
$$
\Ds \vartheta_\a= \calF( \calF(\Ds \vartheta_\a ) )
=\calF( \rho^{2s} \calF(\vartheta_\a)(\rho))= m_\a\calF({\rho^{\frac{-N}{2}+s-\a}})=\g_\a r^{-2s}\vartheta_\a.
$$
For the proof of the fact that the map $\a\mapsto\g_\a$ is continuous and decreasing, we refer to \cite{DDM}.\\
  We define
%
$$
 \Upsilon_\a(t,x) =\begin{cases}
 P(t,\cdot)* \vartheta_\a(t,x) \quad\forall (t,x)\in \R^{N+1}_+\\
\vartheta_\a({x})\quad\forall {x}\in \de \R^{N+1}_+\setminus\{0\},
\end{cases}
$$
where $ P $ is the Poisson  kernel defined in Section \ref{ss:D-to-N}. Clearly $ \Upsilon_\a $ is positive.
 We have that $$-t^{1-2s}\frac{\de \Upsilon_\a}{\de t}=\k_s\Ds \vartheta_\a\quad \textrm{ in }\R^{N}\setminus\{0\}.$$
 Hence we get \eqref{eq:v-al_def}.\\
From  the  regularity theory of \cite{CS},
we deduce that $\Upsilon_\a\in C^\b\left(\ov{\R^{N+1}_+}\setminus\{0\} \right)$ for some $\b>0$. In addition
 $\Upsilon_\a\in H^1\left( \O\times(t_1,t_2); t^{1-2s}\right)$ for every $\O\subset\subset\R^N\setminus\{0\}$ and $0<t_1<t_2<\infty$.\\
Observe that $ \Upsilon_\a(\l z)= \l^{ \frac{2s-N}{2}+\a}  \Upsilon_\a(z)$ and thus choosing $ \l=|z|^{-1}$, we infer that 
 $\Upsilon_\a(z)\leq  \Upsilon_\a(z |z|^{-1})|z|^{ \frac{2s-N}{2}+\a}\leq C |z|^{ \frac{2s-N}{2}+\a}$, for every $z=(t,x)\in B^{N+1}_+(0,R)$ and $R>0$. 
From this, we deduce that  $t^{1-2s}|z|^{-2}\Upsilon_\a^2\in L^1(B^{N+1}_+ (0,R))$ for  $\a>0$.
We also have  $|{x}|^{-2s}\Upsilon_\a^2\in L^1_{loc}(\R^{N})$
for  $\a>0$.\\
We let $\vp$ be a cut-off function such that $\vp=0$ for $|z|<\e$, $\vp=1$
for $2\e<|z|<R$, $\vp=0$ for $|z|>2R$ and $|\n \vp|\leq C \e^{-1}$ for $\e<|z|<2\e $.
 We use $\vp^2 \Upsilon_\a$ as a test function in
\eqref{eq:v-al_def} to get
$$
\int_{\R^{N+1}_+ }t^{1-2s}\n \Upsilon_\a\n(\vp^2 \Upsilon_\a)= \int_{\de \R^{N+1}_+ }|{x}|^{-2s} \Upsilon_\a^2\vp^2.
$$
Integrating by parts  and using Young's inequality, for some constant $c>0$,
we have
$$
c\int_{B^{N+1}_+(0,R)\setminus B^{N+1}_+(0,\e)}t^{1-2s}\vp^2|\n \Upsilon_\a|^2 \leq \int_{B^N(0,2R) }|{x}|^{-2s}\Upsilon_\a^2\vp^2+
 \int_{ \R^{N+1}_+ }t^{1-2s} \Upsilon_\a^2 |\n\vp|^2.
$$
Therefore
$$
c\int_{B^{N+1}_+(0,R)\setminus B^{N+1}_+(0,\e)}t^{1-2s}|\n \Upsilon_\a|^2 \leq\int_{B^N(0,2R)  }|{x}|^{-2s}\Upsilon_\a^2+
 \int_{\e<|(t,x)|<2\e } t^{1-2s}|(t,x)|^{-2}\Upsilon_\a^2.
$$
 Fatou's lemma yields $|\n \Upsilon_\a|\in L^2(B^{N+1}_+(0,R); t^{{1-2s}})$ for $\a>0$.
\QED
%
%
The comparison result  obtained in Lemma \ref{lem:Complem} allows us to derive the following estimate
when the potential $b(x)$ is the Hardy one.
\begin{Lemma}\label{lem:estim}
Let $N>2s$,  $\a\in[0,(N-2s)/2)$ and $p>1$.
 Suppose that $u\in \calL^1_s\cap L^p_{loc}( B^N(0,2)\setminus\{0\}) $, $u\gvertneqq0$
 such that
\be\label{eq:Dsugqg}
\Ds u-  \g_\a |x|^{-2s} u \geq u^p\quad\textrm{ in } \calD'(B^N(0,2)\setminus\{0\} ).
\ee
Then there exists a constant $C>0$ such that
\be\label{eq:coer-ne}
\|\vp\|^2_{\scrH^s_0( B^N(0,2) )}-  \g_\a \int_{ B^N(0,2) }|x|^{-2s}  \vp^2dx \geq
C \int_{ B^N(0,2) }\vp^2dx\quad\forall \vp\in C^{\infty}_c( B^N(0,2) ).
\ee
Moreover there exists a constant $C'>0$ such that
\be\label{eq:ugergst}
u\geq v_\a\geq  C' |x|^{\frac{2s-N}{2}+\a} \quad \textrm{ in }  B^N(0,1),
\ee
where $v_\a\in\scrH^s_{0,b}(B^{N}(0,2)) $ (with $b(x)=\g_\a |x|^{-2s}$) is  the solution to
$$
\calB_s v_\a-\g_\a |x|^{-2s}v_\a=\min(u^p,1)\quad\textrm{ in } B^N(0,2).
$$
\end{Lemma}
\proof
Inequality \eqref{eq:coer-ne} is trivial for $\a>0$ so we consider only the case $\a=0.$
Since $C_s(\{0\})=0$ provided $N>2s$, by Lemma \ref{lem:min-solH10-loc} and Lemma \ref{lem:smmp}, we have
$u>0$ in $B^N(0,2) $ and
$$
M=ess\inf_{B^{N}(0,1) } u>0 .
$$
Hence  $u\in \calL^1_s$ satisfies
$$
\Ds u-  \g_\a |x|^{-2s} u \geq M^{p-1} u\quad\textrm{ in } \calD'(B^N(0,1)\setminus\{0\} ).
$$
By Lemma \ref{lem:AP} we obtain \eqref{eq:coer-ne} thanks to the scale invariance of the  integrals on
the left hand side.\\
We put $ b(x)=\g_\a |x|^{-2s} $. By \eqref{eq:coer-ne},
for $\a\geq0$, we can let $v_\a\in\scrH^s_{0,b}(B^{N}(0,2)) $  be the solution to
$$
\calB_s v_\a-\g_\a |x|^{-2s}v_\a=\min(u^p,1)\quad\textrm{ in } B^N(0,2).
$$
Then by Lemma \ref{lem:Complem},
we have $u\geq v_\a$ in $ B^N(0,2) $.
We first consider the case $\a>0$. Then $\g_\a<\g_0$ and thus $v_\a\in\scrH^s_{0}(B^{N}(0,2)) $.
By the regularity result of \cite{CS}, we get that $V_\a=\calH(\ti{v_\a})$ is continuous
in $\ov{B^{N+1}_+(0,1)}\setminus \{0\}$ and $V_\a\geq0$ by Lemma \ref{lem:mmp1}. Consider $\vartheta_\a$ and its
harmonic extension $\Upsilon_\a$ given by Lemma \ref{lem:grstate}.
  The maximum principal (see Lemma \ref{lem:smmp}) implies that  we can set
\be\label{eq:Ceqmom}
C'= \frac{\displaystyle \min_{\ov{S^N_+}} v_\a  }{\displaystyle   \max_{\ov{S^N_+}} \Upsilon_\a}>0.
\ee
 Put $w=C' \Upsilon_\a-V_\a$. We have weakly
 \be\label{eq:PBgleq}
 \begin{cases}
 \div(t^{1-2s}\n {w})=0\qquad \textrm{ in } B^{N+1}_+(0,1),\\
{w}\leq0 \qquad  \textrm{ on }\ov{S^N_+}, \\
-t^{1-2s}\frac{\de {w}}{\de t}-\k_s\g_\a|{x}|^{-2s}w\leq 0\qquad \textrm{ on } B^N(0,1).
%
\end{cases}
 \ee
Then  $w^+:=\max(w,0)\in H^1_{0,S}(B^{N+1}_+(0,1);t^{1-2s})$
and therefore by integration by parts
$$
\int_{\R^{N+1}_+}t^{1-2s}|\n w^+|^2 dxdt - \k_s\g_\a\int_{B^N(0,1) }|{x}|^{-2s}(w^+)^2dx \leq0.
$$
In particular
$$
\|w^+\|_{\scrH^s_0( B^N(0,1))}^2 - \g_\a\int_{B^N(0,1) }|{x}|^{-2s}(w^+)^2dx \leq0.
$$
Hence $w^+\equiv0$ by Hardy's inequality. Hence  $v_\a\geq C'\vt_\a$ in $ B^N(0,1)$
that is \eqref{eq:ugergst} for $\a>0$.\\
For the case $\a=0$, we put $\a_n=1/n$ and  we notice that the sequence $v_{\a_n}\in \scrH^s_0(B^N(0,2))$
solution to the problem
$$
\calB_s v_{\a_n}-\g_{\a_n} |{x}|^{-2s}v_{\a_n}=f\quad\textrm{ in } B^N(0,2)
$$
is monotone increasing to $v_{0}$ because the mapping $\a\mapsto \g_\a$ is decreasing.
Therefore, taking into account \eqref{eq:Ceqmom}, we readily get \eqref{eq:ugergst}.
\QED

\noindent
\textbf{ Proof of Theorem \ref{thm:main-th-ne}}
\begin{Lemma}\label{lem:ne}
 Let $E$ be a bounded Lipschitz domain of $\R^N$, $N>2s$.
Suppose that $0\in E$ and $\a\in[0,(N-2s)/2)$. Let $u\in \calL^1_s\cap L^p_{loc}(E \setminus\{0\})$ such that 
\begin{equation}
\label{eq:inequalityb-sup}
(-\Delta)^s u-\g_\a |x|^{-2s}\, u \geq  u^p\quad\textrm{ in } \calD'(E\setminus\{0\}),
\end{equation}
with $$\g_\a=2^{2s} \frac{\Gamma\left(\frac{N+2s+2\a}{4}\right)}{\Gamma\left(\frac{N-2s-2\a}{4}\right)} 
\frac{\Gamma\left(\frac{N+2s-2\a}{4}\right)}{\Gamma\left(\frac{N-2s+2\a}{4}\right)} .$$
If $p\geq \frac{N+2s-2\a}{N-2s-2\a}$, then $u=0$ in $E$.
\end{Lemma}
\proof
Assume that $u\neq0$.
It follows from   Lemma \ref{lem:estim}  that there exist  $r, C_r>0$ such that
\be\label{eq:estim-contr}
u(x)\geq v_\a(x) \geq C_r|x|^{\frac{2s-N}{2}+\a}\quad \forall x\in B^N(0,r)\subset E,
\ee
where $v_\a\in\scrH^s_{0,b}(E) $ (with $b(x)=\g_\a |x|^{-2s}$) is  the solution to
$$
\calB_s v_\a-\g_\a |x|^{-2s}v_\a=\min(u^p,1)\quad\textrm{ in } E.
$$
On the other hand Lemma \ref{lem:AP} yields,
for all $\vp\in C^{\infty}_c(E),$

\be\label{eq:APapl}
\|\vp\|^2_{\scrH^s_0(E)}- \g_\a \int_{E}|x|^{-2s}\vp^2dx\geq\int_{E}u^{p-1}\vp^2dx
.
\ee
We first consider the case $\a>0$. If $r$ is small, by \eqref{eq:estim-contr} we have, for $0<\a'<\a$,
 $$
(-\Delta)^s u-\g_{\a'} |x|^{-2s}\, u \geq\left( -\g_{\a'}+\g_\a+ C_r^{p-1}   \right)|x|^{-2s}u\quad\textrm{ in } \calD'(B^N(0,r)\setminus\{0\}).
$$
By Lemma \ref{lem:Complem} and using the same arguments as in Lemma \ref{lem:estim} we get, provided $\a' \nearrow\a$,
\be\label{eq:estim-contr-alpos}
u(x)\geq   C_r'|x|^{\frac{2s-N}{2}+\a'}\quad \forall x\in B^N(0,r/2),
\ee
for some constant $C_r'>0$. Using the   estimate \eqref{eq:estim-contr-alpos}  in \eqref{eq:APapl} we get
$$
\|\vp\|^2_{\scrH^s_0(B^N(0,r/2))}- \g_\a \int_{B^N(0,r/2)}|x|^{-2s}\vp^2dx\geq
(C_r')^{p-1}\int_{B^N(0,r/2)}|x|^{\left(\frac{2s-N}{2}+\a'\right)(p-1)}\vp^2dx,
$$
for all $\vp\in C^{\infty}_c(B^N(0,r/2))$. Since $p> \frac{N+2s-2\a'}{N-2s-2\a'}$, we have
 $$ 
 -\d:=\left(\frac{2s-N}{2}+\a'\right)(p-1)+2s<0.
 $$
 Hence for every $\rho\in(0,r/2)$
$$
\|\vp\|^2_{\scrH^s_0(B^N(0,\rho))}\geq (\g_\a+(C_r')^{p-1} \rho^{-\d} ) \int_{B^N(0,\rho)}|x|^{-2s}\vp^2dx\quad \forall
\vp\in C^{\infty}_c(B^N(0,\rho)).
$$
 This contradicts the sharpness of the  Hardy constant thanks to the scale invariance  of the inequality.\\
Finally, for the case $\a=0$ we note that \eqref{eq:APapl} implies, by density, that
$$
\|v_\a \|^2_{\scrH^s_0(E)}- \g_\a \int_{E}|x|^{-2s}v_\a^2dx\geq\int_{E}v_\a^{p+1}dx.
$$
This  also leads to a contradiction because $v_\a\in \scrH^s_{0,b}(E)$ while by \eqref{eq:estim-contr}
\begin{eqnarray*}
\int_{E}v_\a^{p+1}dx&\geq& C' \int_{B^N(0,r)}|x|^{\left(\frac{2s-N}{2}\right)(p+1)}dx\\
&\geq &C'\int_{S^{N-1}}\int_0^rt^{-1}dtd\s=+\infty,
\end{eqnarray*}
for some constant $C'>0$.
\QED
%
%
\section{Existence of positive solutions}\label{s:exist}
The proof will be separated into several cases. We put
$$
E_\a(u):=\|u\|_{\scrH^s_0(B)}^2-\g_\a\int_{B}|x|^{-2s}u^2dx\quad \forall u\in C^\infty_c(B),
$$
where $B$ is a ball in $\R^N$ centered at 0 with $N>2s$.\\
\textbf{Case 1: $\a\in(0,(N-2s)/2]$ and $1<p<(N+s)/(N-2s)$. }\\
Thanks to the Hardy-Littlewood-Sobolev inequality and Hardy's inequality we have
\be
E_\a(u)\geq \left(\int_{B}u^{p+1}\right)^{2/(p+1)}\quad
\forall u\in C^\infty_c(B).
\ee
Thanks to the compact embedding of $\scrH^s_0(B)$ into $L^{p+1}(B)$, we can minimize
$E_\a$ over the set
\be
\left\{ u\in \scrH^s_{0}(B) \,:\, \int_{B}(u^+)^{p+1}=1\right\}.
\ee
Let  $ u\in \scrH^s_{0}(B)$ be the minimizer. Put $ u^\pm=\max(\pm u,0)$. By Proposition \ref{prop:H1zT-eqscrH},  $u^\pm$  belongs to
$\scrH^s_{0}(B)$. We check rapidly that $E_\a(u^+)\leq E_\a(u) $.  Observe that
\begin{eqnarray*}
\k_s \|u^\pm\|_{\scrH^s_{0}(B)}^2&=&\int_{\R^{N+1}_+}t^{1-2s}|\n \calH(u^\pm)|^2dxdt \leq \int_{\R^{N+1}_+}t^{1-2s}|\n (\calH(u)^\pm)|^2dxdt 
\end{eqnarray*}
because  $ \calH(u^\pm) $ and $\calH(u)^\pm $ have the same trace on $\R^N$ while
$\calH(u^\pm)$ has minimal Dirichlet  energy. Now using this and Hardy's inequality   we have
\begin{eqnarray*}
\k_s  E_\a(u^+)&=&\k_s\|u^+\|_{\scrH^s_{0}(B)}^2-\k_s\g_\a\int_{B}|x|^{-2s}(u^+)^2dx \nonumber \\
&= &\k_s\|u^+\|_{\scrH^s_{0}(B)}^2-\k_s\g_\a\int_{B}|x|^{-2s}u^2dx +\k_s\g_\a\int_{B}|x|^{-2s}(u^-)^2dx\\
&\leq &\k_s\|u^+\|_{\scrH^s_{0}(B)}^2-\k_s\g_\a\int_{B}|x|^{-2s}u^2dx +\k_s\|u^-\|_{\scrH^s_{0}(B)}^2\\
& \leq&\int_{\R^{N+1}_+}t^{1-2s}|\n( \calH(u)^+)|^2dxdt-\k_s\g_\a\int_{B}|x|^{-2s}u^2dx\\
&&+\int_{\R^{N+1}_+}t^{1-2s}|\n( \calH(u)^-)|^2dxdt \label{eq:minHu+}\\
&= &\int_{\R^{N+1}_+}t^{1-2s}|\n \calH(u)|^2dxdt-\k_s\g_\a\int_{B}|x|^{-2s}u^2dx \nonumber \\
&=&\k_s E_\a(u).        
\end{eqnarray*}
Thus we may assume that $u=u^+$ is a nonegative and nontrivial minimizer therefore there  exists a Lagrange multiplier $\l>0$
such that
$$
\calB_s u -\g_\a|x|^{-2s} u=\l u^{p}\quad\textrm{ in } B.
$$
Hence $\l^{\frac{1}{p-1}}\ti{u}$ is a solution of problem \eqref{eq:inequalityb}.\\
\textbf{Case 2: $\a=0$ and $1<p<(N+2s)/(N-2s)$. }\\
{Lemma} \ref{lem:coerciv} yields for every   $ q\in\left(2,\max\left(1,\frac{2}{1+2s}\right)\right) $
$$
 E_0 (u) \geq   \|u\|^{2}_{W^{\t,q}_0(B )} \quad u\in C^\infty_c(B),
$$
with  $ \t=\frac{1+2s}{2}-\frac{1}{q}$. Therefore $\scrH^s_{0,b}(B)$ is compactly
embedded into $L^{p+1}(B)$, with $b=\g_0|x|^{-2s}$. Hence we can minimize $E_0$
 over the set
\be\label{eq:setSchrzb-minimiz}
\left\{ u\in \scrH^s_{0,b}(B) \,:\, \int_{B}(u^+)^{p+1}=1\right\}.
\ee
 We have to check again that $E_0 (u^+)\leq E_0 (u) $. But this can be done by density
and using similar arguments as above. We skip the details. We get a positive  minimizer $u= u^+$ of $E_0$
in the set \eqref{eq:setSchrzb-minimiz}.
We conclude that  $\l^{\frac{1}{p-1}}\ti{u}$ is a  solution to   \eqref{eq:inequalityb}  for some Lagrange multiplier $\l>0$.\\
\textbf{Case 3: $\a\in(0,(N-2s)/2)$ and $(N+2s)/(N-2s)\leq p<{(N+2s-2\a)}/{(N-2s-2\a)}$. }\\
Consider $\vartheta_\b=r^{\frac{2s-N}{2}+\b}$ given by Lemma \ref{lem:grstate} which satisfies
$$
\Ds \vartheta_\b=\g_\b |{x}|^{-2s}\,\vartheta_\b\quad \textrm{ in } { \R^N\setminus\{0\}}.
$$
We look for a solution of the form $w=\mu r^{\frac{-2s}{p-1}}$ with a constant $\mu>0$ to be determined in a minute.
Assume that we can take  $\b\geq 0$ such that $ r^{\frac{2s-N}{2}+\b}=r^{\frac{-2s}{p-1}}$ then
\begin{eqnarray*}
\Ds w&=&\g_{\b} |{x}|^{-2s}w+w^p-w^{p-1}w\\
&=&\g_\a |{x}|^{-2s}w+w^p+(\g_{\b}-\g_\a-\mu^{p-1})|x|^{-2s}w.
\end{eqnarray*}
Since $\b\mapsto \g_\b$ is decreasing, we can choose $\mu^{p-1}=\g_{\b}-\g_\a>0$ provided $\a>\b$.
But note that $\a>\b$  as soon as  $p<(N+2s-2\a)/(N-2s-2\a)$ and $p\geq (N+2s)/(N-2s)$ implies $\b\geq 0$.
In conclusion we have, in $\R^N\setminus \{0\}$,
$$
\Ds w-\g_\a|{x}|^{-2s} w= w^p
$$
and $w\in \calL^1_s\cap L^p(B)$ is a  solution to   \eqref{eq:inequalityb}.
\QED

\section{Appendix}\label{a:remainder}
\subsection{Remainder term for the fractional Hardy inequality}\label{ss:Rt}
Let $E$ be a bounded open set of $\R^N$, $N>2s$, with $0\in E$.
The following (local) Hardy inequality  is a consequence  of \eqref{eq:fracHardy}
\be\label{eq:fracHardyE}
\g_0\int_{E}u^2|x|^{-2s}dx\leq \|u\|^2_{\scrH^s_0(E)} \quad \forall u\in C^\infty_c(E).
\ee
In addition the constant $\g_0$ is optimal.
 Our objective, in this section, is  to improve inequality \eqref{eq:fracHardy} in bounded domains of $\R^N$.\\
Many deal of  work  has been done in improving  the classical Hardy inequality starting  from the work of Brezis-V\'{a}zquez \cite{BV}.
We also quote  \cite{BFT}, \cite{VZ}, \cite{GGM} for related improvements. \\
We shall prove a V\'azquez-Zuazua-type (see  \cite{VZ})  improvement  for the fractional Hardy inequality \eqref{eq:fracHardy}. That is
for    $ 2>q>\max\left(1,\frac{2}{2-a}\right) $,
 there exists a constant $C(E)>0$ such that for all $u\in C^\infty_c(E)$,
$$
C(E) \|u\|^{2}_{W^{\t,q}_0(E )}\leq\|u\|_{\scrH^s_0(E)}^2- \g_0 \int_{E}|x|^{a-1}u^2dx,
$$
where $a=1-2s$ and $ \t=\frac{2-a}{2}-\frac{1}{q}$. The proof requires several preliminary lemmata.\\

Consider the function $\Upsilon_0$ defined in Lemma \ref{lem:grstate} satisfying
\be\label{eq:v-ze_def}
\begin{cases}
\div(t^{1-2s} \Upsilon_0)=0\quad \textrm{ in }\R^{N+1}_+,\\
\Upsilon_0=|{x}|^{\frac{2s-N}{2}} \quad \textrm{ on } \R^{N}\setminus\{0\},\\
-t^{1-2s}\frac{\de \Upsilon_0 }{\de t}=\k_s\g_0 |{x}|^{-2s}\,\Upsilon_0\quad\textrm{ on }\R^{N}\setminus\{0\}.
\end{cases}
\ee
We have seen in the proof of  Lemma \ref{lem:grstate} that
\be\label{eq:vleqz}
|\Upsilon_0(z)|\leq C |z|^{(2s-N)/2}\quad \forall z=(x,t)\in B^{N+1}_+.
\ee
By scale invariance, we have that
\be
\Upsilon_0(z)=R^{(N-2s)/2} \Upsilon_0(R z)\quad\forall R>0.
\ee
This implies the estimate
\be\label{eq:vgeqz}
|\Upsilon_0(z)|\geq C |z|^{(2s-N)/2}\quad \forall z\in B^{N+1}_+
\ee
and also
\be\label{eq:nvleqz}
|\n \Upsilon_0(z)|\leq C |z|^{(2s-N)/2-1}\quad \forall z\in B^{N+1}_+.
\ee
We now prove the following result which  were proved in  \cite{DDM} when $s=1/2$.
\begin{Lemma}\label{eq:imprhlfB}
For every $q\in (1,2)$ there exists a constant $C>0$ such that for all $\vp\in C^\infty_c(B^{N+1})$
\be
C \left(\int_{B^{N+1}_+}t^{\frac{qa}{2}}|\n \vp|^q dz \right)^{2/q}\leq \int_{B^{N+1}_+}t^a|\n \vp|^2 dz- \k_s\g_0 \int_{B^N}|x|^{a-1}\vp^2dx,
\ee
where $a=1-2s$.
\end{Lemma}
\proof
Let $\vp\in C^\infty_c(B^{N+1}\setminus\{0\}) $ and put $\psi=\frac{\vp}{\Upsilon_0}$. Simple computations yield
$$
|\n \vp|^2=|\Upsilon_0\n\psi|^2+\n \Upsilon_0\cdot \n(\Upsilon_0\psi^2).
$$
Integration by parts and using \eqref{eq:v-ze_def} leads to
$$
\int_{B^{N+1}_+}t^a|\n \vp|^2 dz-\k_s\g_0 \int_{B^N}|x|^{a-1}\vp^2dx \geq \int_{B^{N+1}_+}t^a \Upsilon_0^2|\n \psi|^2 dz.
$$
By \eqref{eq:vgeqz} and using polar coordinates $z=r\s=|z|\frac{z}{|z|}$,  we get
\be
\int_{B^{N+1}_+}t^a \Upsilon_0^2|\n \psi|^2 dz\geq C\int_0^1\int_{S^N_+}(\s_1)^a r |\n\psi|^2d\s dr,
\ee
where $\s_1$ is the  component of $\s$ in the $t$ direction.
We wish to show that there exists a constant $C>0$ such that
\be\label{eq:objectif}
I:=\left(\int_0^1\int_{S^N_+}(\s_1)^a r |\n\psi|^2d\s dr\right)^{q/2}\geq C \int_{B^{N+1}_+}t^{\frac{qa}{2}}|\n \vp|^q dz
\quad \forall \vp\in C^\infty_c(B^{N+1}).
\ee
We have
\begin{eqnarray*}
 \int_{B^{N+1}_+}t^{\frac{qa}{2}}|\n \vp|^q dz= \int_{B^{N+1}_+}t^{\frac{qa}{2}}|\n\psi \Upsilon_0+\psi\n \Upsilon_0|^qdz
\leq C \int_{B^{N+1}_+}t^{\frac{qa}{2}}\left(|\n\psi \Upsilon_0|^q+|\psi\n \Upsilon_0|^q\right)dz.
\end{eqnarray*}
Put
$$
I_1=\int_{B^{N+1}_+}t^{\frac{qa}{2}}|\n\psi \Upsilon_0|^qdz,
$$
$$
I_2= \int_{B^{N+1}_+}t^{\frac{qa}{2}}|\psi\n \Upsilon_0|^qdz.
$$
Using \eqref{eq:vleqz} and H\"{o}lder inequality, we have
\begin{eqnarray}\label{eq:compI1}
 I_1 & \leq & C\int_{S^N_+}(\s_1)^{qa/2} \int_0^1 r^{N+q(1-N)/2}|\n \psi|^q drd\s \nonumber\\
&=& C\int_0^1 r^{N+q(1-N)/2-q/2} \int_{S^N_+}(\s_1)^{qa/2} r^{q/2}|\n \psi|^q d\s dr \\
& \leq & C\int_0^1 r^{N(2-q)/2}\left( \int_{S^N_+} (\s_1)^{a} r|\n \psi|^2 d\s   \right)^{q/2}dr\nonumber\\
& \leq & C\int_0^1 \left( \int_{S^N_+} (\s_1)^{a} r|\n \psi|^2 d\s   \right)^{q/2}dr\nonumber\\
& \leq & C \left( \int_0^1\int_{S^N_+} (\s_1)^{a} r|\n \psi|^2 d\s  dr \right)^{q/2}\nonumber\\
& \leq & C I\nonumber.
\end{eqnarray}
On the other hand we have using \eqref{eq:nvleqz}
\begin{eqnarray}\label{eq:compI2}
  I_2\leq &C&\int_{S^N_+}(\s_1)^{qa/2} \int_0^1 r^{N(2-q)/2-q/2}|\psi|^q drd\s\nonumber\\
 &\leq& C \int_{S^N_+}(\s_1)^{qa/2} \int_0^1 r^{N(2-q)/2-q/2+q}\left|\frac{\de \psi}{\de r}\right|^q drd\s \nonumber\\
&\leq& C \int_0^1 r^{N+q(1-N)/2-q/2} \int_{S^N_+}(\s_1)^{qa/2} r^{q/2}|\n \psi|^q d\s dr\\
& \leq & C I\nonumber,
\end{eqnarray}
where in the second inequality we have used the one dimensional
Hardy inequality $\int_0^1 f^q dr\leq c \int_0^1 r^{-q}|f'|^q dr$ and observing that
\eqref{eq:compI1} is just \eqref{eq:compI2}. The lemma follows because  $C^\infty_c(B^{N+1}\setminus\{0\})$
is dense in $C^\infty_c(B^{N+1}) $  with respect to the $H^{1}(B^{N+1} ;t^{a})$-norm when $N\geq 3$, see \cite{Kil}.
\QED
%
The Lion's interpolation inequality, [\cite{Lions} Paragraph 5], shows that
for $a\neq 0$ and $-\frac{1}{q}<\frac{a}{2}<\frac{1}{q}$ there exits a constant $C>0$ such that
\be\label{eq:JLL-inter}
C \|v\|^q_{W^{\t,q}(\R^N)}\leq \int_{\R^{N+1}_+ }t^{\frac{a q}{2}}|\n v|^qdxdt
+\int_{\R^{N+1}_+ }t^{\frac{a q}{2}}| v|^qdxdt,\quad\forall v\in C^\infty_c(\R^{N+1}),
\ee
where $ \t=\frac{2-a}{2}-\frac{1}{q}$.\\
We first prove a generalized weighted Poincar\'e trace inequality. The proof is standard.
We recall that the space $  W^{1,q}_{0,S}(B^{N+1}_+;t^{qa/2 })$ was defined in Section \ref{s:NP}.
\begin{Lemma}\label{eq:Gen_poinct-trace}
Suppose that $q>1$ and   $\t:=\frac{2-a}{2}-\frac{1}{q}>0$. Then the following inequality holds
\be\label{eq:GP}
\int_{B^{N+1}_+}t^{qa/2}|u|^qdz\leq C\int_{B^{N+1}_+}t^{qa/2}|\n u|^qdz+ C \int_{B^{N}}|u|^qdz\quad\forall  u\in W^{1,q}_{0,S}(B^{N+1}_+;t^{qa/2 }).
\ee
\end{Lemma}
\proof
Inequality \eqref{eq:GP}  is well known for $a=0$. So we  restrict ourself to the case $a\neq0$.\\
Assume by contradiction that \eqref{eq:GP} does not hold.  Then there exits a sequence  $u_n\in W^{1,q}_{0,S}(B^{N+1}_+;t^{qa/2 }) $ such that
\be\label{eq:nrj_un-to0}
\int_{B^{N+1}_+}t^{qa/2}|\n u_n|^qdz+  \int_{B^{N}}|u_n|^qdz=o(1)
\ee
and
$$
\int_{B^{N+1}_+}t^{qa/2}|u_n|^qdz>0\quad \forall n\in\N.
$$
Up to  normalization, we may assume that $ \int_{B^{N+1}_+}t^{qa/2}|u_n|^qdz=1$. But then \eqref{eq:nrj_un-to0}
implies that $u_n$ is bounded
in $ W^{1,q}(B^{N+1}_+;t^{qa/2 })$ thus $u_n \rightharpoonup u$ in $W^{1,q}(B^{N+1}_+;t^{qa/2 })$
 and $u_n\to u$ in $L^q(B^{N+1}_+ ;t^{qa/2 })$ (see \cite{Froh}) so that
\be\label{eq:nm_ueq1}
 \int_{B^{N+1}_+}t^{qa/2}|u|^qdz=1.
\ee
It follows from \eqref{eq:JLL-inter}  and the compact embedding of $ W^{\t,q}_{0}(B^{N}) $ into  $ L^q(B^N)$ ,
 that  $ u_n \to u$ in $ L^q(B^N)$. From \eqref{eq:nrj_un-to0} we  get  $u\Big|_{\R^N}=0$ and also $\n u=0$.
 It turns out that $u=0$ in $ B^{N+1}_+$, a contradiction with \eqref{eq:nm_ueq1}.
\QED
By an argument of partition of unity, we have  the following
\begin{Lemma}\label{eq:imprhlfBL2}
Let  $ 2>q>\max\left(1,\frac{2}{2-a}\right) $. There exist some constants $C,c>0$ such that for all $\vp\in C^\infty_c(\R^{N+1})$
\be
C \|\vp\|^{2}_{W^{\t,q}(B^N )}\leq\int_{\R^{N+1}_+}t^a|\n \vp|^2 dz- \k_s\g_0 \int_{\R^N}|x|^{a-1}\vp^2dx + c\int_{\R^N}\vp^2dx,
\ee
where $a=1-2s$ and $  \t=\frac{2-a}{2}-\frac{1}{q}$.
\end{Lemma}
\proof
We put
$$
J(v)=\int_{\R^{N+1}_+}t^a|\n v|^2 dz- \k_s\g_0 \int_{\R^N}|x|^{a-1}v^2dx.
$$
Let $\chi\in C^\infty_c (B^{N+1})$, $0\leq\chi\leq 1$ in $B^{N+1}$ and  such that $\chi\equiv1$ on $B^{N+1}(0,1/2)$.
 Let $\eta\in H^1(B_+^{N+1};t^{a})$ be the minimum of the problem
 $$\inf\left\{\int_{B^{N+1}_+}t^{a}|\n u|^2dz\,:\, u-\chi \in H^1_0(B^{N+1}_+; t^a)\right\}.$$ Then
$$
\begin{cases}
 \div(t^{a}\n \eta)=0\quad\textrm{ in } B_+^{N+1}\\
\eta=\chi \quad\textrm{ on } B^N\\
\eta=0 \quad\textrm{ on } \ov{S^{N}_+}.
\end{cases}
$$
It turns out that $0\leq \eta\leq 1$ in $B_+^{N+1}$. In addition, thanks to \cite{CSS}, $\lim_{t\to 0} t^a\frac{\de \eta}{\de t}\in L^\infty_{loc}(B^N)$.
Given $\vp\in C^\infty_c(\R^{N+1})$, simple computations based on integration by parts lead to
$$
\int_{\R^{N+1}_+}t^a|\n(\vp\eta)|^2dz\leq \int_{\R^{N+1}_+}t^a|\n \vp|^2dz +c \int_{B^N}\vp^2dx,
$$
where $c>0$ depends only on $\eta$. On the other hand we have
\begin{eqnarray*}
\int_{B^N}|x|^{a-1}(\eta\vp)^2dx&=&\int_{\R^N}|x|^{a-1}\vp^2dx +  \int_{\R^N}(1-\eta^2)|x|^{a-1}\vp^2\\
&\leq &\int_{ \R^N}|x|^{a-1}\vp^2dx + c\int_{\R^N}\vp^2dx.
\end{eqnarray*}
Therefore we obtain
$$
J(\vp\eta) \leq J(\vp)+   c\int_{\R^N}\vp^2dx.
$$
Applying Lemma \ref{eq:imprhlfB}, we infer that
$$
\left(\int_{\R_+^{N+1}}t^{\frac{qa}{2}}|\n (\eta \vp)|^q dz \right)^{2/q}\leq  J(\vp)+ c \int_{\R^N}\vp^2dx .
$$
By Lemma \ref{eq:Gen_poinct-trace} and H\"{o}lder inequality ($1<q<2$)
$$
C \|\eta \vp\|^{2}_{W^{1,q}(\R^{N+1}_+;t^{qa/2})}\leq \ J(\vp)
 \displaystyle+ c \int_{\R^N}\vp^2dx.
$$
Using Lions' interpolation inequality \eqref{eq:JLL-inter} with $a\neq0$, we obtain
$$
C \|\eta\vp\|^{2}_{W^{\t,q}(\R^N )}\leq  J(\vp) \displaystyle+ c \int_{\R^N}\vp^2dx.
$$
If $a=0$, it is well know that $W^{1,q}(\R^{N+1}_+) $ embeds continuously into $ W^{1-1/q,q}(\R^N )$.
Recalling that $\eta=\chi\equiv 1$ on $B^N(0,1/2)$, the lemma follows by scaling.
\QED
Taking advantages to the singular nature of the Hardy potential and the scale invariance, we prove the main result in this section:
\begin{Lemma} \label{lem:coerciv}
Let  $ 2>q>\max\left(1,\frac{2}{2-a}\right) $.
Then there exists a constant $C_0>0$ such that for all $u\in C^\infty_c(B^{N})$,
\be
C_0 \|u\|^{2}_{W^{\t,q}_0(B^N )}\leq\|u\|_{\scrH^s_0(B^N)}^2- \g_0 \int_{B^N}|x|^{a-1}u^2dx,
\ee
where $a=1-2s$ and $ \t=\frac{2-a}{2}-\frac{1}{q}$.
\end{Lemma}
\proof
Let $u\in C^\infty_c(B^{N})$ and we define $U=\calH(\ti{u})= \calH({u})$. Then $U\in H^1(\R^{N+1}_+;t^a)$  thus by Lemma \ref{eq:imprhlfBL2} and a
 density argument we get
$$
C \|U \|^{2}_{W^{\t,q}(B^N )}\leq\int_{\R^{N+1}_+}t^a|\n U|^2 dz- \k_s\g_0 \int_{\R^N}|x|^{a-1}U^2dx+c  \int_{\R^N}U^2dx.
$$
Since $ u=U$ on $B^N$, it follows that
$$
C \|u\|^{2}_{W^{\t,q}(B^N )}\leq\|u\|_{\scrH^s_0(B^N)}^2- \g_0 \int_{B^N}|x|^{a-1}u^2dx+c  \int_{B^N}u^2dx.
$$
Let $r\in(0,1)$ we derive from the above that for every $u\in C^\infty_c(B^{N}(0,r)) $
\be\label{eq:almst-impr}
C \|u\|^{2}_{W^{\t,q}(B^N(0,r) )}\leq\|u\|_{\scrH^s_0(B^N(0,r))}^2- \g_0 \int_{B^N(0,r)}|x|^{a-1}u^2dx+c  \int_{B^N(0,r)}u^2dx,
\ee
with $c,C>0$ independent on $r$.
This holds because $\|u\|_{\scrH^s_0(B^N(0,r))}=\|u\|_{\scrH^s_0(B^N(0,1))} $ as long as $u\in C^\infty_c(B^{N}(0,r)) $ and $r<1$.\\

It is clear from \eqref{eq:almst-impr} that wee need only to
 show that  there exists a constant $C_2>0$ such that for every  $u\in C^\infty_c(B^{N}) $
\be\label{Improv_obkec-Poinc}
C_2 \int_{B^N}u^2dx\leq \|u\|_{\scrH^s_0(B^N)}^2- \g_0 \int_{B^N}|x|^{a-1}u^2dx.
\ee
Before proceeding, we recall that the mapping $\a\mapsto \g_\a$ (defined in Lemma \ref{lem:grstate}) is decreasing.
 We will use this fact and the estimates in Lemma  \ref{lem:estim} to conclude the proof.
Pick
$$
\a\in \left(0,\frac{N-2s}{2} \right)
$$
and let $r>0$ be so small  that
\be\label{eq:gz-ga*r-cpssit}
(\g_0-\g_\a) r^{a-1}-c>0.
\ee
As we did in Section \ref{s:exist},
by \eqref{eq:almst-impr}, we can  define the space $\scrH^s_{0,b}(B^N(0,r))$ with $b(x)= \g_0|x|^{a-1}-c$.
Letting $2<p+1<\frac{2N}{N-2s}$, we can choose $q$ (close to 2) so that  $W^{\t,q}_0(B^N(0,r) )$
is compactly embedded into  $ L^{p+1}(B^N(0,r) ) $. Then minimization procedure implies that
there exits a nonnegative and nontrivial  $u_r\in \scrH^s_{0,b}(B^N(0,r))$ solution to
$$
\calB_s u_r+c u_r-\g_0|x|^{a-1}u_r = C u_r^p \quad\textrm{ in } B^N(0,r).
$$
By density
\be \label{u_rsolv}
\Ds\ti{ u_r}+c \ti{ u_r}-\g_0|x|^{a-1}\ti{u_r} = C \ti{u_r}^p\quad\textrm{ in } \calD'(B^N(0,r)).
\ee
We have, by \eqref{eq:gz-ga*r-cpssit},
$$
\Ds \ti{u_r}-\g_\a|x|^{a-1}\ti{u_r}\geq ((\g_0-\g_\a) r^{a-1}-c) \ti{u_r}+ C \ti{u_r}^p\geq  C\ti{ u_r}^p \quad\textrm{ in }\calD'( B^N(0,r)).
$$
Therefore by Lemma  \ref{lem:estim},  there exists a constant $C_r>0$ such that
\be\label{eq:estm-u_r}
u_r\geq C_r |x|^{-\frac{N-2s}{2}+\a }\quad\textrm{ in } B^N(0,r/2).
\ee
For every $k\in \N$, take $v^k\in\scrH^s_{0}(B^N(0,r)) $ as the   solution to
$$
\calB_s v_k+c v_k-(\g_0-1/k)|x|^{a-1}v_k= C \min( u_r^p,k)\quad\textrm{ in } B^N(0,r).
$$
 Since $\ti{u_r}$ satisfies \eqref{u_rsolv}, it follows (see Remark \ref{rem:complem-c}) that
\be\label{eq:u_rgev_k}
u_r\geq v_k>0 \quad \textrm{ in }B^N(0,r)
\ee
by Lemma \ref{lem:mmp1} and Lemma \ref{lem:smmp}.\\
Put $V_k=\calH(\ti{v_k})$ we have that for any $\Psi\in H^1_{0,T}( B^N(0,r); t^{a}) $
\be\label{eq:VkPsi}
\begin{array}{ccc}
\k_s^{-1}\displaystyle\int_{\R^{N+1}_+}t^a\n V_k \cdot\n \Psi dz&= &(\g_0 -1/k)\displaystyle\int_{B^N(0,r)}|x|^{a-1}v_k \Psi dx
\displaystyle -c  \int_{B^N(0,r)}v_k \Psi dx\\
&&+\displaystyle C \int_{B^N(0,r)}\min(u_r^p,k) \Psi dx.
\end{array}
\ee
Thanks to \eqref{eq:estm-u_r}, we can choose $r'\in(0,r/2)$ (small) such that
\be\label{eq:rprchose}
-c+C u_r^{p-1}\geq  1 \quad\textrm{ in } B^N(0,r').
\ee
For such a fixed $r'$, take $\vp\in C^\infty_c( B^N(0,r'))$. For $\e>0$, set $V_k^\e= V_k+\e$
and put  $\psi:=\frac{\calH(\vp)}{V_k^\e} $.
We have
$$
|\n \calH(\vp)|^2=|V_k^\e\n\psi|^2+\n V_k\cdot \n({V_k^\e}\psi^2)
$$
and also $V_k^\e\psi^2\in  H^1_{0,T}( B^N(0,r); t^{a}) $. This together with \eqref{eq:VkPsi} yields
%
\begin{eqnarray*}
\k_s^{-1} \int_{\R^{N+1}_+ }t^{a}|\n \calH(\vp) |^2dxdt&\geq &
(\g_0 -1/k)\int_{B^N(0,r')}|x|^{a-1}\frac{v_k}{v_k+\e}\vp^2dx\\
&& -c  \int_{B^N(0,r')}\frac{v_k}{v_k+\e} \vp^2 dx\\
& & +C \int_{B^N(0,r')}\frac{\min(u_r^p,k)}{v_k+\e} \vp^2 dx.
\end{eqnarray*}
Thus
\begin{eqnarray*}
\k_s^{-1} \int_{\R^{N+1}_+ }t^{a}|\n \calH(\vp) |^2dxdt&\geq &
(\g_0 -1/k)\int_{B^N(0,r')}|x|^{a-1}\frac{v_k}{v_k+\e}\vp^2dx -c  \int_{B^N(0,r')} \vp^2 dx\\
& & +C \int_{B^N(0,r')}\frac{\min(u_r^p,k)}{v_k+\e} \vp^2 dx.
\end{eqnarray*}
By Fatou's lemma, when $\e\to 0$, we have
\begin{eqnarray*}
\k_s^{-1} \int_{\R^{N+1}_+ }t^{a}|\n \calH(\vp) |^2dxdt&\geq &
(\g_0 -1/k)\int_{B^N(0,r')}|x|^{a-1}\vp^2dx-c  \int_{B^N(0,r')} \vp^2 dx\\
& & +C \int_{B^N(0,r')}\frac{\min(u_r^p,k)}{v_k} \vp^2 dx.
\end{eqnarray*}
By \eqref{eq:u_rgev_k} we obtain
\begin{eqnarray*}
\k_s^{-1} \int_{\R^{N+1}_+ }t^{a}|\n \calH(\vp) |^2dxdt&\geq &(\g_0 -1/k)\int_{B^N(0,r')}|x|^{a-1}\vp^2dx -c  \int_{B^N(0,r')} \vp^2 dx\\
& & +C \int_{B^N(0,r')}\frac{\min(u_r^{p},k)}{u_r} \vp^2 dx.
\end{eqnarray*}
It follows again from Fatou's lemma that
\begin{eqnarray*}
 \k_s^{-1}\int_{\R^{N+1}_+ }t^{a}|\n \calH(\vp) |^2dxdt&\geq &\g_0 \int_{B^N(0,r')}|x|^{a-1}\vp^2dx +  \int_{B^N(0,r')}(-c+ C u_r^{p-1}) \vp^2 dx.\\
\end{eqnarray*}
From the choice of $r'$ in \eqref{eq:rprchose}, we get immediately, for any  $\vp\in C^\infty_c( B^N(0,r'))$,
$$
\k_s^{-1} \int_{\R^{N+1}_+ }t^{a}|\n \calH(\vp) |^2dxdt-\g_0 \int_{B^N(0,r')}|x|^{a-1}\vp^2dx\geq \int_{B^N(0,r')} \vp^2 dx.
$$
By scaling we have  \eqref{Improv_obkec-Poinc} which was our objective.

\QED


   \textbf{ Acknowledgments} \\
   This work is supported by the Alexander-von-Humboldt Foundation. The author would like to thank Professor Tobias Weth
for fruitful discussions.

    \label{References}


\begin{thebibliography}
   \footnotesize
   %

%
   \bibitem{AS-CPDE} S. Armstrong, B. Sirakov,   Nonexistence of positive supersolutions of elliptic equations via the maximum principle.
   Comm. PDE, to appear.
   \bibitem{AS-RIMS} S. Armstrong, B. Sirakov, A new approach to Liouville theorems for elliptic inequalities. RIMS proceeding, to appear.
   %
   \bibitem{BG} P. Baras and J. A. Goldstein, The heat equation with a singular potential,
   Trans. Amer. Math. Soc.  284 (1984), no. 1, 121-139.
   %
%
%
%
\bibitem{BFT} G. Barbatis, S. Filippas, A.  Tertikas, A unified approach to
improved $L^p$ Hardy inequalities with best constants . Trans. Amer.
Math. Soc., 356, (2004), 2169-2196.
%
%
\bibitem{BBC}K. Bogdan; K. Burdzy; Z-Q. Chen, Censored stable processes.  Probab. Theory Related Fields  127  (2003),  no. 1, 89-152.
%
   \bibitem{BC} H. Brezis and X. Cabr\'e. Some simple nonlinear PDEs without solutions. Bull. UMI 1 (1998),
   223-262.
   %
   \bibitem{BDT}H. Brezis, L. Dupaigne, A. Tesei,
   On a semilinear elliptic equation with inverse-square potential.
   Selecta Math. (N.S.) 11 (2005), no. 1, 1-7.
   %
   \bibitem{BM} H. Brezis and M. Marcus, Hardy's inequalities revisited.
    Dedicated to Ennio De Giorgi.
     Ann. Scuola Norm. Sup. Pisa Cl. Sci. (4)  25  (1997),  no. 1-2, 217-237.
   %
    \bibitem{BMS} H. Brezis, M. Marcus and I. Shafrir,
      Extermal functions for Hardy's inequality with weight, J. Funct. Anal. 171 (2000), 177-191.
   %
\bibitem{BV} H. Brezis and J. L. V\'azquez, Blow-up solutions of some nonlinear elliptic problems, Rev.
    Mat. Univ. Complut. Madrid 10 (1997), no. 2, 443-469.

\bibitem{CS} X. Cabr\'e and Y. Sire, Nonlinear equations for fractional Laplacians I: Regularity, maximum
    principles, and Hamiltonian estimates (2010).

%
   \bibitem{CT} X. Cabr\'e and  J. Tan,  Positive solutions of nonlinear
   problems involving the square root of the Laplacian.  Adv. Math.  224  (2010),  no. 5, 2052-2093.
   %
\bibitem{CSilv} L. Caffarelli and L. Silvestre, An extension problem related to the fractional Laplacian, Comm.
     Partial Differential Equations 32 (2007), no. 7-9, 1245-1260.
\bibitem{CSS} L. A. Caffarelli, S. Salsa, and L. Silvestre, Regularity estimates for the solution and the free
     boundary of the obstacle problem for the fractional Laplacian, Invent. Math. 171 (2008),
     no. 2, 425-461.

   \bibitem{CDDS} A. Capella, J. D\'{a}vila, L. Dupaigne, Y. Sire.  Regularity of radial extremal solutions for some non local semilinear equations. Preprint 2010.
   http://www.capde.cl/publication/.
   %
   \bibitem{DD-H}J. D\'{a}vila; L. Dupaigne,  Hardy-type inequalities.
    J. Eur. Math. Soc. (JEMS) 6 (2004), no. 3, 335-365.
   %
   %
\bibitem{DDM}] J. D\'avila, L. Dupaigne, and M. Montenegro, The extremal solution of a boundary reaction
             problem, Commun. Pure Appl. Anal. 7 (2008), no. 4, 795-817.
%


   %
   %

%
   \bibitem{D}L. Dupaigne, Semilinear elliptic PDEs with a singular potential, J. Anal. Math., 86
   (2002), 359-398.
   %
   %
   \bibitem{DN} L. Dupaigne and G. Nedev. Semilinear elliptic PDEs with a singular potential. Adv. Differential Equations 7 (2002), 973-1002.
   %
\bibitem{FKS} E. B. Fabes, C. E. Kenig, and R. P. Serapioni, The local regularity of solutions of degenerate
     elliptic equations, Comm. Partial Differential Equations 7 (1982), no. 1, 77-116.

%
   \bibitem{mmf}M. M. Fall, On the Hardy Poincar\'e inequality with boundary singularities.
   Commun. Contemp. Math, to appear.

   %
   \bibitem{FaMu1}M. M. Fall, R. Musina,  Hardy-Poincar\'e inequality with boundary singularities.
   Proc. Roy. Soc. Edinburgh, to appear.
   %
   \bibitem{FaMu-ne} M. M. Fall, R. Musina,  Sharp nonexistence results for a linear elliptic
    inequality involving Hardy and Leray potentials.
   J. Inequal.  Appl. 2011 (2011). doi:10.1155/2011/917201.
   %
\bibitem{Fall-ne-sl}M. M. Fall,  Nonexistence of distributional supersolutions of a semilinear elliptic equation with Hardy potential (preprint 2011).
%
\bibitem{FQT} P. Felmer, A. Quaas, J. Tan, Positive solutions of Nonlinear Schr\"{o}dinger equation with the fractional Laplacian. Preprint 2010.

%
\bibitem{Froh} A. Fr\"{o}hlich, The Helmholtz decomposition of weighted $L^q$-spaces for Muckenhoupt weights.
  Ann. Univ. Ferrara Sez. V. 46, 1 (2000), 11-19.
%
\bibitem{GGM} F. Gazzola, H.-C. Grunau, E. Mitidieri, Hardy inequalities with optimal constants and remainder terms,
Trans. Amer. Math. Soc. 356 (2004), 2149-2168.

\bibitem{Gris} P. Grisvard,  Elliptic problems in nonsmooth domains. Monographs and Studies in Mathematics, 24. Pitman (Advanced Publishing Program), Boston, MA, 1985.
  %
\bibitem{Herb} I. W. Herbst,  Spectral theory of the operator $(p^{2}+m^{2})^{1/2}-Ze^{2}/r$.  Comm. Math. Phys.  53  (1977), no. 3, 285-294.
%
\bibitem{Kil} T. Kilpel\"{a}inen,
Weighted Sobolev spaces and capacity.
Ann. Acad. Sci. Fenn. Ser. A I Math. 19 (1994), no. 1, 95-113.
%
   \bibitem{KLS-Trans}V. Kondratieva, V. Liskevich  Z.  Sobolb, Positive supersolutions to semi-linear second-order non-divergence type elliptic equations in exterior domains.
   Trans. Amer. Math. Soc. Volume 361 (2009), 697-713.
   %
   \bibitem{KLS}V. Kondratieva, V. Liskevich  Z.  Sobolb,
   Positive solutions to semi-linear and quasi-linear second-order elliptic equations on unbounded domains. Handbook of  Differential  Equation.
   Stationary Partial Differential Equations,
   volume 6 (Edited by M. Chipot) 2008 Elsevier, pp.177-268.
   %
\bibitem{Landk} N. S. Landkof. Foundations of modern potential theory. Springer-Verlag,
New York, 1972. Translated from the Russian by A. P. Doohovskoy, Die
Grundlehren der mathematischen Wissenschaften, Band 180.

%
\bibitem{Lions}  J. L. Lions,  Th\'eor\'{e}mes de trace et d'interpolation. I.  Ann. Scuola Norm. Sup. Pisa (3)  13  1959 389-403.
%
   \bibitem{LLM} V. Liskevich, S. Lyakhova and V. Moroz, Positive solutions to singular semilinear elliptic
   equations with critical potential on cone-like domains, Adv. Differential Equations 11 (2006), pp. 361-398.
   %
   %
   \bibitem{LLM-Proce} V. Liskevich, S. Lyakhova and V. Moroz, Positive solutions to semilinear elliptic equations with critical lower order terms.
   EQUADIFF 2003, 549-554, World Sci. Publ., Hackensack, NJ, 2005.
   %
   \bibitem{KLS-JDE}V. Kondratieva, V. Liskevich  Z.  Sobolb, Second-order semilinear elliptic inequalities in exterior domains.
   J. of Differential Equations, vol. 187 (2003), 429-455.
   %
\bibitem{Maz}V.G. Maz'ja, Sobolev Spaces, Springer-Verlag, 1985.

   %
   \bibitem{PoTe}S. I. Pohozaev and A. Tesei, Nonexistence of local solutions to semilinear partial differential
   inequalities. Ann. Inst. H. Poincar\'e Anal. Non Lin\'eaire 21 (2004), pp. 487-502.
   %
   %
\bibitem{Rup} R. L. Frank, A Simple Proof of Hardy-Lieb-Thirring Inequalities. Commun. Math. Phys. 290, 789-800 (2009).
%
   \bibitem{Sil}L. Silvestre, Regularity of the obstacle problem for a fractional power of the
   Laplace operator. Comm. Pure Appl. Math., 60(1):67-112, 2007.
%
\bibitem{SV} Y. Sire, E. Valdinoci, Fractional Laplacian phase transitions and boundary reactions: a geometric inequality and a symmetry result,
  Journal of Functional Analysis, 256, 6 (2009), 1842-1864.
%
\bibitem{StWe} E. M. Stein and G. Weiss,  Fractional integrals on $n$-dimensional Euclidean space.  J. Math. Mech.  7  1958 503-514.
%
\bibitem{SW} E. M. Stein  and G. Weiss,
Introduction to Fourier analysis on Euclidean spaces.
Princeton Mathematical Series, No. 32. Princeton University Press, Princeton, N.J., 1971.
   %
%
   \bibitem{T} S. Terracini, On positive entire solutions to a class of equations with a singular coefficient and critical exponent.
    Adv. Differential Equations  1  (1996),  no. 2, 241-264.
   %
\bibitem{VZ}J. L. V\'azquez, E. Zuazua, The Hardy inequality and the asymptotic behavior of the heat equation with an inverse-square potential,
 J. Funct. Anal. 173 (2000).

%

%
\bibitem{Yaf}D. Yafaev, Sharp constants in the Hardy-Rellich inequalities.
J. Funct. Anal. 168 (1999), no. 1, 121-144.
   \end{thebibliography}
   \end{document}